\documentclass[11pt,a4paper]{amsart}
\usepackage{amssymb,amsmath,epsfig,graphics,mathrsfs}
\usepackage{mathtools}
\usepackage[backend=biber]{biblatex}
\usepackage{amssymb,amsthm,epsfig,esint}
\usepackage{fancyhdr}
\usepackage{hyperref}
\hypersetup{colorlinks=true, urlcolor=black, linkcolor=black, citecolor=red, bookmarksopen=true}
\usepackage{amsmath}
\usepackage{amsfonts}
\usepackage{amssymb}
\usepackage{amsthm}
\usepackage{epsfig,graphics,mathrsfs}
\usepackage{graphicx}

\usepackage[usenames, dvipsnames]{color} 

\usepackage{hyperref}

\theoremstyle{plain}
\newtheorem{thrm}{Theorem}[section]
\newtheorem*{thrm*}{Theorem}
\newtheorem{lemma}[thrm]{Lemma}
\newtheorem{prop}[thrm]{Proposition}
\newtheorem{cor}[thrm]{Corollary}
\theoremstyle{definition}
\newtheorem{dfn}[thrm]{Definition}
\theoremstyle{remark}

\newtheoremstyle{named}{}{}{\normalfont}{}{\bfseries}{}{ }{}
\theoremstyle{named}
\newtheorem*{namedthm}{}

\numberwithin{equation}{section}
\usepackage{color}

\begin{document}

\newcommand{\fdelta}{D_if_{\delta}}
\newcommand{\Dfdelta}[1]{{D}_{#1}f_\delta}
\newcommand{\Dfdeltaij}{D^2_{i,j}f_\delta}
\newcommand{\Dfdeltainl}{D^2_{i,\ell+s_\ell n}f_\delta}
\newcommand{\Dfdeltanl}{D_{\ell+s_\ell n}f_\delta}
\newcommand{\Dfdeltajnl}{D^2_{j,\ell+s_\ell n}f_\delta}

\newcommand{\tx}{\tilde x}
\newcommand{\R}{\mathbb R}
\newcommand{\N}{\mathbb N}
\newcommand{\C}{\mathbb C}
\newcommand{\lie}{\mathcal G}
\newcommand{\hN}{\mathcal N}
\newcommand{\D}{\mathcal D}
\newcommand{\A}{\mathcal A}
\newcommand{\B}{\mathcal B}
\newcommand{\sL}{\mathcal L}
\newcommand{\sLi}{\mathcal L_{\infty}}

\newcommand{\G}{\Gamma}
\newcommand{\x}{\xi}

\newcommand{\eps}{\epsilon}
\newcommand{\al}{\alpha}
\newcommand{\be}{\beta}
\newcommand{\p}{\partial}  %followed by _
\newcommand{\lig}{\mathfrak}

\def\dist{\mathop{\varrho}\nolimits}

\newcommand{\BCH}{\operatorname{BCH}\nolimits}
\newcommand{\Lip}{\operatorname{Lip}\nolimits}
\newcommand{\Hol}{C}
\newcommand{\lip}{\operatorname{lip}\nolimits}
\newcommand{\capQ}{\operatorname{Cap}\nolimits_Q}
\newcommand{\pCap}{\operatorname{Cap}\nolimits_p}
\newcommand{\Om}{\Omega}
\newcommand{\om}{\omega}
\newcommand{\half}{\frac{1}{2}}
\newcommand{\e}{\varepsilon}
\newcommand{\vn}{\vec{n}}
\newcommand{\X}{\Xi}
\newcommand{\tLip}{\tilde  Lip}

\newcommand{\Span}{\operatorname{span}}

\newcommand{\ad}{\operatorname{ad}}
\newcommand{\Hm}{\mathbb H^m}
\newcommand{\Hn}{\mathbb H^n}
\newcommand{\Hone}{\mathbb H^1}
\newcommand{\Lie}{\mathfrak}
\newcommand{\Layer}{V}
\newcommand{\hgrad}{\nabla_{\!H}}
\newcommand{\im}{\textbf{i}}
\newcommand{\nz}{\nabla_0}
\newcommand{\s}{\sigma}
\newcommand{\se}{\sigma_\e}

\newcommand{\ued}{u^{\e,\delta}}
\newcommand{\ueds}{u^{\e,\delta,\sigma}}
\newcommand{\tnabla}{\tilde{\nabla}}

\newcommand{\bx}{\bar x}
\newcommand{\by}{\bar y}
\newcommand{\bt}{\bar t}
\newcommand{\bs}{\bar s}
\newcommand{\bz}{\bar z}
\newcommand{\btau}{\bar \tau}

\newcommand{\LC}{\mbox{\boldmath $\nabla$}}
\newcommand{\Ne}{\mbox{\boldmath $n^\e$}}
\newcommand{\nuo}{\mbox{\boldmath $n^0$}}
\newcommand{\nuu}{\mbox{\boldmath $n^1$}}
\newcommand{\nue}{\mbox{\boldmath $n^\e$}}
\newcommand{\nuek}{\mbox{\boldmath $n^{\e_k}$}}
\newcommand{\dse}{\nabla^{H\Su, \e}}
\newcommand{\dso}{\nabla^{H\Su, 0}}
\newcommand{\tX}{\tilde X}

\newcommand\red{\textcolor{red}}
\newcommand\green{\textcolor{green}}

\newcommand{\Xie}{X^\epsilon_i}
\newcommand{\Xje}{X^\epsilon_j}
\newcommand{\Su}{\mathcal S}
\newcommand{\F}{\mathcal F}

\title[]{Gradient estimates for an orthotropic nonlinear diffusion equation in the Heisenberg group}

\author{Michele Circelli}


\begin{abstract}
We establish that weak solutions to a parabolic orthotropic $p$-Laplacian-type equation in the Heisenberg group $\Hn$ exhibit local Lipschitz continuity with respect to the spatial variables, uniformly in time, for the range $2\leq p\leq4$.
\end{abstract}

\maketitle

\tableofcontents

\section*{Introduction}

In this paper, we extend the results of \cite{circelli2024lipschitz} to the non-stationary setting, adapting the methodologies introduced in \cite{capogna2020regularity} and \cite{capogna2021lipschitz}. Specifically, we consider an open bounded set $\Omega$ in the Heisenberg group $\mathbb{H}^n$ and a time horizon $T > 0$, and we study the quasilinear, degenerate parabolic equation
\begin{equation}\label{maineq0}
    \partial_t u = \textnormal{div}_H ( D f(\nabla_H u) ), \quad \text{in} \quad Q := \Omega \times (0, T),
\end{equation}
where $f: \mathbb{R}^{2n} \to \mathbb{R}$ is defined by
\begin{equation}\label{ortotropicH}
f(z) = \frac{1}{p} \sum_{i=1}^{n} \left( z_i^2 + z_{i+n}^2 \right)^{\frac{p}{2}},
\end{equation}
with $p \geq 1$, and $Df = (D_1 f, \ldots, D_{2n} f)$ denotes the Euclidean gradient of $f$. Here, $\nabla_H u = (X_1 u, \ldots, X_{2n} u)$ represents the horizontal gradient of a weak solution $u \in L^p((0, T), HW^{1,p}(\Omega))$, where $HW^{1,p}(\Omega)$ denotes the Sobolev space associated with $\nabla_H$.

Defining
\begin{equation*}
\aligned
    &\lambda_{i}(z) = \left( z_i^2 + z_{i+n}^2 \right)^{\frac{p-2}{2}}, \quad \forall i \in \{1, \ldots, n\}, \\
    &\lambda_{i}(z) = \left( z_{i-n}^2 + z_i^2 \right)^{\frac{p-2}{2}}, \quad \forall i \in \{n+1, \ldots, 2n\},
\endaligned
\end{equation*}
equation \eqref{maineq0} can be equivalently rewritten as
\begin{equation}\label{maineq}
\partial_t u = \sum_{i=1}^{2n} X_i (\lambda_i (\nabla_H u) X_i u), \quad \text{in} \quad Q = \Omega \times (0, T).
\end{equation}

We prove that weak solutions to \eqref{maineq} are locally Lipschitz continuous with respect to the spatial variables, uniformly in time, for the range $2 \leq p \leq 4$. As byproducts, we also obtain local $L^q$ integrability results for both the non-horizontal derivative and the time derivative of solutions.

Our main result is the following.

\begin{namedthm}
\textbf{Main Theorem.}  
Let $2\leq p\leq4$ and let $u\in L^p((0,T), HW^{1,p}(\Om))$ be a weak solution to \eqref{maineq} in $Q=\Om\times (0,T)$, then
\begin{equation*}
    \nabla_Hu\in L_{loc}^\infty(Q,\mathbb{R}^{2n}).
\end{equation*}
More precisely, for every parabolic cylinder 
\begin{equation*}
    Q_{\mu,r}(x_0,t_0):=B(x_0,r)\times(t_0-\mu \,r,t_0)
\end{equation*}
such that $Q_{\mu,2r}(x_0,t_0)\subset\subset Q$ there is a constant $c=c(n, p,L)>0$ such that
\begin{equation}\label{c1alpha}
    \|\nabla_Hu\|_{L^\infty(Q_{\mu,r}(x_0,t_0),\, \R^{2n})}\leq c\mu^{\frac{1}{2}}\max\left\{
    \bigg(\frac{1}{\mu r^{N+2}}\int\int_{Q_{\mu,2r}(x_0,t_0)}|\nabla_H u|^p\, dx dt \bigg)^{\frac{1}{2}} ,\mu^{\frac{p}{2(2-p)}}\right\},
\end{equation}
where $N:=2n+2$ is the homogeneous dimension of $\mathbb{H}^n$.

In addition $\partial_tu,Zu\in L_{loc}^q(\Omega),$ for any $1\leq q<\infty$. 
\end{namedthm}

The present paper, along with \cite{capogna2020regularity} and \cite{capogna2021lipschitz}, represents the first instances in the literature that study higher regularity for weak solutions of non-stationary degenerate quasilinear equations in the sub-Riemannian setting. In these references, the authors investigate the regularity of solutions to equation \eqref{maineq0}, modeled on the $p$-Laplace equation, i.e. when 
\begin{equation*}
    f(z)=\frac{1}{p}\left(\delta+|\nabla_Hu|^2\right)^{\frac{p}{2}}.
\end{equation*}
In \cite{capogna2020regularity}, the authors prove the smoothness of solutions in the non-degenerate case, namely $\delta>0$, while in \cite{capogna2021lipschitz}, they establish the Lipschitz regularity for solutions in the degenerate case ($\delta=0$), with $2 \leq p \leq 4$. Both studies are based on techniques introduced by Zhong in \cite{zhong2017regularity}, where the optimal regularity for solutions in the stationary case is addressed.

However, none of these results apply to our equation \eqref{maineq0}. Indeed, all of them rely on the fact that the loss of ellipticity of the operator $\textnormal{div}_H Df$ occurs only at a single point $z = 0$, whereas equation \eqref{maineq} is much more degenerate, as it degenerates in the unbounded set
\begin{equation*}
	\bigcup_{i=1}^n\{ \lambda_i (z) = 0 \} = \bigcup_{i=1}^n \left\{ z_i^2 + z_{i+n}^2 = 0 \right\} \subset \mathbb{R}^{2n},
\end{equation*}
which is the union of $2n-2$ dimensional submanifolds. For this reason, we do not expect Hölder regularity of the gradient in the spatial variable, uniformly in time, but only its boundedness. Such result, in the Euclidean setting, was proven in \cite{brascoparabolic}: the Main Theorem in this paper shows that this is also the case in the setting of the Heisenberg group for the range $2 \leq p \leq 4$. We believe that the boundedness of the gradient is true for the entire range $1 < p < \infty$, but it requires different techniques. For a comprehensive overview of degenerate parabolic equations in the Euclidean setting, modeled on the $p$-Laplacian, we refer to \cite{dibenedettolibro}.

The study of the stationary case of equation \eqref{maineq}, motivated by its relation to the optimal transport problem with congestion in the Heisenberg group (see \cite{circelli2023transport} and \cite{circelli2024continuous}), was recently addressed in \cite{circelli2024lipschitz}, where the authors proved the local Lipschitz regularity for solutions for the range $p \geq 2$, adapting to the orthotropic case the techniques introduced by Zhong in \cite{zhong2017regularity}.

Regarding the stationary case of equation \eqref{maineq} in the Euclidean setting, it is worth mentioning \cite{bousquet2018lipschitz}, where the Lipschitz regularity for solutions is addressed for $p \geq 2$, even in a more degenerate case, and \cite{Demengel} for an alternative proof based on viscosity methods. See also \cite{bousquet2020lipschitz} and \cite{bousquet2023singular} for the same result in the anisotropic case. In the particular case of the plane, Bousquet and Brasco \cite{bousquet2018c1} proved that weak solutions are $C^1$ for $1 < p < \infty$. Moreover, derivatives of solutions have a logarithmic modulus of continuity: see \cite{ricciotti2018regularity} for the case $1 < p < 2$, and \cite{lindqvist2018regularity} for the case $p \geq 2$. 

\bigskip

The proof of the Main Theorem builds upon the extension of the results in \cite{circelli2024lipschitz} to the parabolic setting, incorporating a Poincaré-type inequality for the vertical derivative of smooth functions (Lemma \ref{poincarétypeineq0}), introduced in \cite{capogna2021lipschitz} and inspired by techniques originally developed for the study of the Levi mean curvature fully nonlinear PDE in \cite{cittimontanari}. As a first step, we approximate equation \eqref{maineq} via a regularized equation obtained through the Riemannian approximation of the Heisenberg group, for which the smoothness of solutions follows directly from classical regularity theory for parabolic equations.

A central result in establishing the local Lipschitz continuity of solutions is the Caccioppoli-type estimate for the first derivatives of the approximating solutions, stated in Theorem \ref{maincaccioppoli}. This estimate is uniform with respect to the approximating parameters and, when combined with the aforementioned Poincaré-type inequality, yields a uniform local $L^q$ bound for the non-horizontal derivatives of the approximating solutions (Lemma \ref{Poincaretypeineq}). Notably, this is the only point in the analysis where the restriction $2 \leq p \leq 4$ is required.

From Theorem \ref{maincaccioppoli}, a local uniform bound for the gradient of the approximating solutions is obtained via a Moser-type iteration, ultimately leading to the Lipschitz regularity of solutions when passing to the limit, as described in Section \ref{S:proofmainthm}. Furthermore, also the local $L^q$ integrability of the non-horizontal derivatives of solutions, for any $1 \leq q < \infty$, follows by letting the approximating parameters tend to zero. Finally, using a standard argument—independent of the Heisenberg group structure—we also establish the local $L^q$ integrability of the time derivative of solutions, for any $1 \leq q < \infty$. 

\bigskip

The plan of the paper is as follows. Section \ref{S:prelim} is devoted to preliminary results: we introduce the Heisenberg group and its Riemannian approximation, and establish the notations used throughout the paper. In Section \ref{S:approximatingequation}, we approximate the main equation with a more regular one, employing a smooth uniformly convex approximation of the orthotropic function \eqref{ortotropicH} together with the Riemannian approximation of the Heisenberg group. Section \ref{S:cacc} is dedicated to proving Caccioppoli-type inequalities for the derivatives of the approximating solutions, which involve the non-horizontal derivatives.

In Section \ref{S:integrabilityZu}, we derive a uniform (in the approximating parameters) integrability estimate for the vertical derivative of the approximating solutions, using a Poincaré-type inequality. Section \ref{S:maincaccioppoli} builds on these results to establish the main Caccioppoli-type inequality for the first derivatives of the approximating solutions (Theorem \ref{maincaccioppoli}). This inequality is uniform in the approximating parameters and no longer involves the non-horizontal derivatives.

In Section \ref{S:lipschitzestimate}, we apply a Moser iteration scheme, built on the main Caccioppoli-type inequality proved in the previous section, to obtain a local uniform bound on the $L^\infty$ norm of the approximating solutions. Section \ref{S:integtimederivative} focuses on a uniform integrability estimate for the time derivative of the approximating solutions. Finally, in Section \ref{S:proofmainthm}, we prove the main theorem by passing to the limit as the approximating parameters tend to zero, from which the local Lipschitz continuity with respect to the spatial variables, uniformly in time, follows as in the classical Euclidean setting.

\section{Preliminaries}\label{S:prelim}

\label{section:preliminaries}

In this section we fix our notation: we introduce the Heisenebrg group $\mathbb{H}^n$ and we collect some preliminary results that will be used throughout the rest of the paper. 

\subsection{The Heisenberg group}

Let $n\geq 1$, we identify the Heisenberg group $\Hn$ with the Euclidean space $\mathbb{R}^{2n+1}$, equipped with the group multiplication
\begin{equation*}
    xy=\left(x_1+y_1, \dots, x_{2n}+y_{2n}, z+s+\frac{1}{2}\sum_{i=1}^n (x_iy_{i+n}-x_{i+n}y_i)\right),
\end{equation*}
for any two points
$x=(x_1,\ldots,x_{2n},z), y=(y_1,\ldots,y_{2n},s)\in {\mathbb H}^n$. 

The left invariant vector fields corresponding to the
canonical basis of the Lie algebra
\[ X_i=\partial_{x_i}-\frac{x_{i+n}}{2}\partial_t, \quad
X_{i+n}=\partial_{x_{i+n}}+\frac{x_i}{2}\partial_t,\quad 1\leq i\leq n,\] 
are also called horizontal vector fields and we denote by
$\nabla_Hu=\sum_{i=1}^{2n}X_iuX_i\cong \left(X_1u,\ldots,X_{2n}u\right)$ the horizontal gradient of any smooth function $u:\mathbb{H}^n\to\mathbb{R}$. Given a smooth horizontal vector field $\phi = \sum_{i=1}^{2n} \phi_i X_i$, its horizontal divergence is $\textnormal{div}_H \phi = \sum_{i=1}^{2n} X_i \phi_i$. 

The only
non-trivial commutator
\[ Z=\partial_z=[X_i,X_{i+n}]=X_iX_{i+n}-X_{i+n}X_i,\quad 1\le i\le n,\]
is also called vertical (or non-horizontal) vector field. We denote by $N:=2n+2$ the homogeneous dimension of $\mathbb{H}^n$.

Let us introduce the Carnot-Carathéodory distance $d$. An absolutely continuous curve $\gamma:[0,1]\to\mathbb{R}^{2n+1}$ is said to be horizontal if its tangent vector $\dot{\gamma}(t)\in\textnormal{span}\left(X_1(\gamma(t)),\ldots,X_{2n}(\gamma(t))\right)$, at almost every $t\in[0,1]$. Due to the stratification of the space the H\"ormander condition is satisfied, and the  Rashevsky-Chow's theorem guarantees that any couple of points can be joined with an horizontal curve  \cite{Chow}. It is then possible to give the following definition of distance. Given $x,y\in\mathbb{H}^n$, the Carnot-Carathéodory distance between them is defined as
\begin{equation*}\label{CCdistance}
	d(x,y):= \inf \  \left\{  \int_0^1|\dot{\gamma}(t)|\, dt \ : \  \gamma \text{ is horizontal}, \ \gamma(0)=x, \ \gamma(1)=y \right\}, 
\end{equation*}
where $|\cdot|$ denotes the norm associated with the left-invariant sub-Riemannian metric $g$, defined by $g\left(X_i,X_j\right)=\delta_{i,j}$, for $i,j=1,\ldots,2n$.

All of the balls are defined with respect to the Carnot-Carathéodory distance: $$B(x,r)=\{ y\in {\mathbb H}^n: d(x,y)<r\},\, x\in\mathbb{H}^n,r>0.$$

The Haar measure in ${\mathbb H}^n$ is the Lebesgue measure of ${\mathbb R}^{2n+1}$. 

Let $1 \leq p < \infty$ and let $\Omega \subset \mathbb{H}^n$ be an open set. We denote by $HW^{1,p}(\Omega)$ the Sobolev space of functions $u \in L^p(\Omega)$ whose horizontal gradient $\nabla_H u$ belongs to $L^p(\Omega, \mathbb{R}^{2n})$. This space is endowed with the norm
\begin{equation*}
    \|u\|_{HW^{1,p}(\Omega)} = \|u\|_{L^p(\Omega)} + \|\nabla_H u\|_{L^p(\Omega, \mathbb{R}^{2n})}.
\end{equation*}
We denote by $HW^{1,p}_0(\Omega)$ the closure of $C^\infty_c(\Omega)$ in $HW^{1,p}(\Omega)$ with respect to this norm.

\subsection{Riemannian approximation of the Heisenberg group}\label{SS:riemannianapprox}

The left-invariant sub-Riemannian structure of $\mathbb{H}^n$ arises as the pointed Hausdorff-Gromov limit of Riemannian manifolds, in which the non-horizontal direction is increasingly penalized.

Let $\e>0$, we denote by $g_\e$ the left-invariant Riemannian metric for which the frame defined by $$X^\e_1 := X_1,...,X^\e_{2n} := X_{2n} \textnormal{ and } X_{2n+1}^\e := \e Z$$ is orthonormal. 

It has been proved by Gromov in \cite{Gromov} that the left invarian Riemannian manifolds $\left(\mathbb{H}^n,g_\e\right)$ converge to the left invariant sub-Riemannian manifold $\left(\mathbb{H}^n,g\right)$, as $\e\to0^+$, in the pointed Hausdorff-Gromov sense, i.e. the $g_\e$-Riemannian balls $B_\e$  satisfy $B_\e\to B$, as $\e\to0^+$, in the Hausdorff-Gromov sense.

If $u:\mathbb{H}^n\to\mathbb{R}$ is any smooth function, the gradient associated with the Riemannian metric $g_\e$ is
\begin{equation*}
	\nabla_\e u:=\sum_{i=1}^{2n+1}X_i^\e uX_i^\e=\sum_{i=1}^{2n}X_iuX_i+\e^2ZuZ\cong\left(X_1u,\ldots,X_{2n}u,\e Zu\right).
\end{equation*}

For a smooth vector field $\phi= \sum_{i=1}^{2n+1} \phi_i X_i^\e$ will also denote 
$$\textnormal{div}_\e\phi = \sum_{i=1}^{2n+1} X_i^\e \phi_i.$$
Formally we have
\begin{equation*}
    \nabla_\e u\to \nabla_Hu,\quad \textnormal{div}_\e\phi \to \textnormal{div}_H \phi,\quad\text{ as }\e\to0^+.
\end{equation*}
We note explicitly that, again formally, we have
\begin{equation*}
    |\nabla_\e u|:=|\nabla_\e u|_{g_\e}=\sum_{i=1}^{2n}(X_iu)^2+\e^2(Zu)^2\to |\nabla_Hu|,\quad \textnormal{as }\e\to0^+.
\end{equation*}

Let $1 \leq p < \infty$ and let $\Omega \subset \mathbb{H}^n$ be an open set. We adopt the unconventional notation $W^{1,p,\varepsilon}(\Omega)$ to denote the Sobolev space of functions $u \in L^p(\Omega)$ whose Riemannian gradient $\nabla_\varepsilon u$ belongs to $L^p(\Omega, \mathbb{R}^{2n+1})$. This space is endowed with the norm
\begin{equation*}
    \|u\|_{W^{1,p,\varepsilon}(\Omega)} = \|u\|_{L^p(\Omega)} + \|\nabla_\varepsilon u\|_{L^p(\Omega, \mathbb{R}^{2n+1})}.
\end{equation*}
We denote by $W^{1,p,\e}_0(\Omega)$ the closure of $C^\infty_c(\Omega)$ in $W^{1,p,\e}(\Omega)$ with respect to this norm.

For our purposes it is important to recall the following version of the Rellich-Kondrachov theorem in this setting. It is a special case of \cite[Theorem 1.3.1]{Ivanov2011}.
\begin{thrm}\label{Rellichthm}
    Let $1\leq p<N$, then the embedding
    \begin{equation*}
        W^{1,p,\e}(\Omega)\xhookrightarrow{}\xhookrightarrow{}L^r(\Omega)
    \end{equation*}
    is compact for any $q$, with $1\leq r<p^*$.
\end{thrm}

\section{An approximating equation}\label{S:approximatingequation}

The proof of the Main Theorem is based on uniform a priori estimates for solutions to a regularized partial differential equation that approximates \eqref{maineq}. The approximation procedure introduced below is a regularization scheme widely used in the literature to derive a priori estimates for weak solutions of PDEs in the Heisenberg group. It relies heavily on the Riemannian approximation of the Heisenberg group; see, for instance, \cite{capogna2020regularity,capogna2021lipschitz,capogna2023regularity,capogna2019conformality} and \cite{circelli2024lipschitz}.

Let $\Omega \subset \mathbb{H}^n$ be an open bounded set, $T > 0$, we define $Q:=\Omega\times (0,T)$ and we consider the equation
\begin{equation}\label{maineq01}
    \partial_tu=\textnormal{div}_H ( D f(\nabla_H u) ),\quad \quad  \text{ in  \quad }Q,
\end{equation}
where $f:\mathbb{R}^{2n}\to\mathbb{R}$ is the function
\begin{equation*}
f(z)=\frac{1}{p}\sum_{i=1}^{n}\left(z_i^2+z_{i+n}^2\right)^{\frac{p}{2}},
\end{equation*}
$p\geq1$. If we denote by 
\begin{equation*}
\aligned
    &\lambda_{i}(z) = \left(z_i^2+z_{i+n}^2\right)^{\frac{p-2}{2}},\quad \forall i\in\{1,\ldots,n\},\\
    &\lambda_{i}(z) = \left(z_{i-n}^2+z_i^2\right)^{\frac{p-2}{2}},\quad \forall i\in\{n+1,\ldots,2n\},
\endaligned
\end{equation*}
then 
\begin{equation*}
    Df(z)=\Big(\lambda_1(z)z_1,\ldots,\lambda_{2n}(z)z_{2n}\Big)
\end{equation*}
and therefore equation \eqref{maineq0} reads as 
\begin{equation}\label{maineq1}
\partial_tu=\sum_{i=1}^{2n}X_i(\lambda_i(\nabla_H u) X_iu),\quad \quad  \text{ in  \quad }Q=\Omega\times (0,T),
\end{equation}
and it satisfies the following structure condition:
\begin{equation*}
    \sum_{i=1}^{n}\lambda_i(z)\left(\xi_i^2+\xi_{i+n}^2\right)\leq\langle D^2f(z)\xi,\xi\rangle\leq (p-1)\sum_{i=1}^{n} \lambda_i(z)\left(\xi_i^2+\xi_{i+n}^2\right).
\end{equation*}

We say that a function $u\in L^p\left((0,T), HW^{1,p}(\Om)\right)$
is a weak solution of \eqref{maineq01} if
\begin{equation}\label{weak}
    \int_0^T \int_\Omega u\, \p_t\phi\  dx dt=\int_0^T \int_\Omega \sum_{i=1}^{2n} D_if(\nabla_H u) X_i \phi \ dx dt,
\end{equation}
for every $\phi\in C^{\infty}_0(Q)$.

As in \cite{circelli2024lipschitz}, we denote by $f_\delta:\mathbb{R}^{2n+1}\to\mathbb{R}$ the function
$$f_{\delta}(z) = \frac{1}{p}\sum_{i=1}^{n} \left(\delta+z_i^2+z_{i+n}^2+ z_{2n+1}^2\right)^{\frac{p}{2}},$$
with $\delta>0$, and we consider the Riemannian parabolic equation 
\begin{equation}\label{approx1}
    \p_tu=\textnormal{div}_\e( D f_\delta(\nabla_\e u) ),
\end{equation}
where $Df_\delta$ denotes the Euclidean gradient of $f_\delta$ in $\mathbb{R}^{2n+1}$ and $\textnormal{div}_\e$ and $\nabla_\e$ are the Riemannian divergence and gradient, respectively, see Subsection \ref{SS:riemannianapprox}.
 
If we denote by 
\begin{equation}\label{lambdaideltaeps}
\begin{split}
    &\lambda_{i,\delta}(z):=(\delta+z_i^2 + z_{i+n}^2+z_{2n+1}^2)^{(p-2)/2},\quad \forall i\in\{1,\ldots,n\},\\
    &\lambda_{i,\delta}(z):=(\delta+z_{i-n}^2 + z_{i}^2+z_{2n+1}^2)^{(p-2)/2},\quad \forall i\in\{n+1,\ldots,2n\},\\&\lambda_{2n+1, \delta}(z):= \sum_{i=1}^{n}\lambda_{i, \delta}(z),
\end{split}
\end{equation}
then the gradient of $f_\delta$ can be written as follows
\begin{equation}\label{structurecondition02}
    Df_{\delta}(z) = \Big(\lambda_{1,\delta}(z) z_1, \cdots, \lambda_{2n,\delta}(z) z_{2n},  \lambda_{2n+1,\delta}(z) z_{2n+1}\Big).
\end{equation}

Hence, the explicit expression of the regularized equation 
becomes 
$$\p_tu = \sum_{i=1}^{2n+1}X_i^\e(\lambda_{i, \delta}(\nabla_\e u) X_i^\e u), 
$$
and the structure condition
\begin{equation}\label{structurecondition01}
	\sum_{i=1}^{n} \lambda_{i,\delta}(z)\left(\xi_i^2+\xi_{i+n}^2+\xi_{2n+1}^2\right)\leq\langle D^2f_{\delta}(z)\xi,\xi\rangle \leq L\sum_{i=1}^{n} \lambda_{i,\delta}(z)\left(\xi_i^2+\xi_{i+n}^2+\xi_{2n+1}^2\right),    
\end{equation}
for any $\xi\in\mathbb{R}^{2n+1}$, where $L=L(n,p)>1$ is a constant.

The next lemma presents a pointwise inequality that will be instrumental in the proof of the Main Theorem. The idea behind the proof originates from \cite[Lemma 4.4]{Marcellini1991}.

\begin{lemma}\label{8agosto2025}
There is a constant $c=c(p,L)>0$ such that, for every $z,\zeta\in\mathbb{R}^{2n+1}$,
\begin{equation*}
    |z_i|^p\leq c\left(\left(\delta+|\zeta|^2\right)^{\frac{p}{2}}+\sum_{i=1}^{2n+1}D_if_\delta(\xi)(z_i-\zeta_i)\right),\quad \forall i=1,\ldots, 2n+1.
\end{equation*}
\end{lemma}

\begin{proof}
For any $z,\zeta\in\mathbb{R}^{2n+1}$ we denote by
\begin{equation*}
    F(t)=\sum_{i=1}^{2n+1}D_if_\delta(tz+(1-t)\zeta)(z_i-\zeta_i).
\end{equation*}
The structure condition \eqref{structurecondition01}, the definition of the $\lambda_{i,\delta}$'s \eqref{lambdaideltaeps} and the Jensen inequality imply
\begin{equation*}
    \aligned
    \sum_{i=1}^{2n+1}\left(D_if_\delta(z)-D_if_\delta(\zeta)\right)(z_i-\zeta_i)&=F(1)-F(0)=\int_0^1\frac{d}{dt}\left(F(t)\right)\ dt\\
    &=\int_0^1\sum_{i,j=1}^{2n+1}\Dfdeltaij(tz+(1-t)\zeta)(z_i-\zeta_i)(z_j-\zeta_j)\ dt\\&\geq\sum_{i=1}^{2n+1} (z_i-\zeta_i)^2\int_0^1\lambda_{i,\delta}(\zeta+t(z-\zeta))\ dt\\
    &\geq\sum_{i=1}^{2n+1}(z_i-\zeta_i)^2\int_0^1\left(\delta+(\zeta_i+t(z_i-\zeta_i))^2\right)^{\frac{p-2}{2}}\ dt\\&\geq\sum_{i=1}^{2n+1}(z_i-\zeta_i)^2\left(\delta+\left(\int_0^1 (\zeta_i+t(z_i-\zeta_i))\ dt\right)^2\right)^{\frac{p-2}{2}}\\
    &=\sum_{i=1}^{2n+1}(z_i-\zeta_i)^2\left(\delta+\left(\frac{z_i+\zeta_i}{2}\right)^2\right)^{\frac{p-2}{2}}.
    \endaligned
\end{equation*}
Using the following inequality, whose proof is contained in the proof from \cite[Lemma 4.4]{Marcellini1991},
\begin{equation*}
    |z_i|^p\leq c\left((z_i-\zeta_i)^2|\zeta_i+z_i|^{p-2}+|\zeta_i|^p\right),\quad \forall i=1,\ldots,2n+1,
\end{equation*}
with $c=c(p)>0$, we may conclude that there is $c_1=c_1(p)>0$ such that 
\begin{equation}\label{11agosto2025}
    |z_i|^p\leq c_1(p)\left(|\zeta|^p+\sum_{i=1}^{2n+1}\left(D_if_\delta(z)-D_if_\delta(\zeta)\right)(z_i-\zeta_i)\right).
\end{equation}

Moreover, if we define the function
\begin{equation*}
    G_i(t):=D_if_\delta(t\zeta),
\end{equation*}
for any $i=1,\ldots,2n+1$, $\zeta\in\R^{2n+1}$, then
\begin{equation}\label{11agosto20251}
    \aligned
    |D_if_\delta(\zeta)|\leq \underbrace{|D_if_\delta(0)|}_{=0}&+\int_0^1|G_i'(t)|dt=\int_0^1\left|\sum_{j=1}^{2n+1}\Dfdeltaij(t\zeta)\zeta_j\right|\ dt\\&\leq L|\zeta_i|^2\int_0^1\lambda_{i,\delta}(t\zeta)\ dt\leq L(\delta+|\zeta|^2)^{\frac{p-1}{2}}.
    \endaligned
\end{equation}
Therefore, by \eqref{11agosto2025} and \eqref{11agosto20251}, we can conclude that there is $c_2=c_2(p,L)>0$ such that, for any $\tau>0$, 
\begin{equation}
    |z_i|^p\leq c_2\left(|\zeta|^p+\sum_{i=1}^{2n+1}D_if_\delta(z)(z_i-\zeta_i)+\frac{\tau^p}{p}\left(|z_i|^p+|\zeta_i|^p\right)+\frac{p-1}{p\tau^{\frac{p}{p-1}}}\left(\delta+|\zeta|^2\right)^{\frac{p}{2}}\right).
\end{equation}
Choosing $\tau>0$ small enough the thesis follows.
\end{proof}

\bigskip

\begin{dfn}[$\delta,\e$-weak solution]
We say that a function $u_{\delta,\e}\in L^p\left((0,T), W^{1,p,\e}(\Omega)\right)$ is a weak solution to the equation \eqref{approx1} if 
\begin{equation}\label{weakformulation}
    \int_0^T \int_\Omega u_{\delta,\e}\, \p_t\phi\  dx dt=\int_0^T\int_{\Omega}\sum_{i=1}^{2n+1}D_i f_{\delta}(\nabla_\e u_{\delta,\e})X_i^\e\phi\ dxdt,
\end{equation}
for every $\phi\in C^\infty_c(Q)$.
\end{dfn}

Since \eqref{approx1} is strongly parabolic (in the Riemannian sense) for every $\delta, \varepsilon > 0$, classical regularity theory for Riemannian parabolic equations (see, for example, \cite{uraltsevaparabolic}) implies that the approximating solutions, namely the weak solutions to \eqref{approx1}, denoted by $u_{\delta,\varepsilon}$, are smooth in $Q$, i.e., $u_{\delta,\varepsilon} \in C^\infty(Q)$. These solutions converge, in a suitable sense, to a function $u_0$ that solves the equation \eqref{maineq01}. With this in mind, the goal of the remainder of the paper is to establish a uniform Lipschitz bound for the approximating solutions, independent of the parameters $\varepsilon$ and $\delta$, so that this estimate can be transferred to solutions of \eqref{maineq01}.

In the next lemma, we collect the partial differential equations satisfied by the first-order derivatives of the approximating solutions $X_\ell^\varepsilon u_{\delta,\varepsilon}$, for $\ell \in \{1, \ldots, 2n+1\}$.

\begin{lemma} 
Let $u_{\delta,\e}$ be a weak solution to \eqref{approx1} in $Q$ and let us denote by $v_\ell = X^\e_\ell u_{\delta,\e}$, with $\ell = 1,\ldots, 2n+1$. 

If $\ell\in\{1,\ldots, 2n\}$, then the function $v_\ell$ solves the equation
\begin{equation}\label{eqX}
    \p_t v_\ell=\sum_{i,j=1}^{2n+1} X^\e_i\left( \Dfdeltaij (\nabla_\e u_{\delta,\e})X^\e_\ell X^\e_ju_{\delta,\e} \right)+s_\ell Z (\Dfdeltanl(\nabla_\e u_{\delta,\e})),
\end{equation}
where $s_\ell = (-1)^{[\frac{\ell}{n+1}]}$ and $[\cdot]$ denotes the floor function;

If $\ell=2n+1$, then $v_{2n+1}$ solves
\begin{equation}\label{eqX2n+1}
    \partial_tv_{2n+1}=\sum_{i,j=1}^{2n+1}X^\e_i(\Dfdeltaij(\nabla_\e u_{\delta,\e})X^\e_\ell X^\e_ju_{\delta,\e}).
\end{equation}
\end{lemma}

For the proof see \cite[Lemma]{capogna2020regularity}.

\section{A Caccioppoli-type inequality for the first derivatives of approximating solutions}\label{S:cacc}

The aim of the next three sections is to derive higher regularity estimates for the approximating solutions $u_{\delta,\varepsilon}$ that remain stable as $\varepsilon, \delta \to 0^+$. Throughout these sections, with a slight abuse of notation, we will omit the indices $\varepsilon, \delta$ and denote by $u$ a weak solution to \eqref{approx1}. Furthermore, we denote by $c$ a positive constant, which may vary from line to line. Unless explicitly stated otherwise, $c$ depends only on the dimension $n$, the exponent $p$, and the constant $L$ appearing in the structure condition \eqref{structurecondition01}. Importantly, we emphasize once again that $c$ does not depend on the approximating parameters $\varepsilon$ and $\delta$, and thus it does not degenerate as $\varepsilon, \delta \to 0^+$.

We begin by establishing a uniform Caccioppoli-type estimate for the first-order derivatives of the approximating solutions $u$, which also involves the vertical derivative $Z u$. The term containing $Z u$ will be removed in Theorem \ref{maincaccioppoli}.

\begin{lemma} \label{CaccioppoliXu} 
Let $u$ be a weak solution to \eqref{approx1} in $Q$. There is $c=c(n,p,L) > 0$ such that, for any $\beta \geq 0$ and any non-negative function $\eta\in C^{1}([0,T], C^\infty_c (\Om))$ vanishing on the parabolic boundary of $Q$, one has
\begin{equation}\label{CaccioppoliXu1}
    \aligned
    &\frac{1}{\beta+2}\sup_{0<t<T}\int_\Om\eta^2 (\delta+|\nabla_\e u|^2)^{\frac{\beta+2}{2}}\ dx\\&+\int_0^T\int_\Om\eta^2 (\delta+|\nabla_\e u|^2)^{\frac{\beta}{2}}\sum_{i=1}^{n}\lambda_{i,\delta}(\nabla_\e u)\left(|\nabla_\e X^\e_i u|^2+|\nabla_\e X^\e_{i+n} u|^2+|\nabla_\e X^\e_{2n+1} u|^2\right)\ dxdt\\
    &\leq c(\beta + 1)\int_0^T\int_\Om \left(|\nabla_\e\eta|^2+\eta|Z\eta|\right)(\delta+|\nabla_\e u|^2)^{\frac{\beta+2}{2}}\sum_{i=1}^n\lambda_{i,\delta}(\nabla_\e u)\ dxdt\\
    &+\frac{2}{\beta+2}\int_0^T\int_\Om\eta|\p_t\eta|(\delta+|\nabla_\e u|^2)^{\frac{\beta+2}{2}}\ dx\\&+c(\beta +1)^2\int_0^T\int_\Om \eta^2 (\delta+|\nabla_\e u|^2)^{\frac{\beta}{2}} \sum_{i=1}^n\lambda_{i,\delta}(\nabla_\e u)|Z u|^2\ dxdt.
    \endaligned
\end{equation}
\end{lemma}

\begin{proof} 

We fix $\eta\in C^{1}([0,T], C^\infty_c (\Om))$ and we want to use $\phi = \eta^2\left(\delta+|\nabla_\e u|^2\right)^{\frac{\beta}{2}}X^\e_\ell u$, with $\ell\in\{1,\ldots,2n+1\}$ as test function in the weak formulation of the equations solved by the first derivatives, i.e. the weak formulation of \eqref{eqX} and \eqref{eqX2n+1}.

If $\ell\in\{1,\ldots,2n\}$, we use $\phi$ as test test function in \eqref{eqX} and, through an integration by parts in the right hand side, we obtain
\begin{align*}
    &\frac{1}{2} \int_0^T\int_\Om \eta^2\left(\delta+|\nabla_\e u|^2\right)^{\frac{\beta}{2}}\p_t \left[ (X^\e_\ell u)^2\right]\ dxdt\\ &+\sum_{i,j=1}^{2n+1}  \int_0^T\int_\Om \Dfdeltaij(\nabla_\e u) X^\e_\ell X^\e_j u  X^\e_i\left(\eta^2\left(\delta+|\nabla_\e u|^2\right)^{\frac{\beta}{2}}X^\e_\ell u\right)\ dxdt\\&=-s_{\ell} \int_0^T\int_\Om \Dfdeltanl(\nabla_\e u)Z\left(\eta^2\left(\delta+|\nabla_\e u|^2\right)^{\frac{\beta}{2}} X^\e_\ell u\right)\ dxdt. 
\end{align*}
We compute $X^\e_i\left(\eta^2\left(\delta+|\nabla_\e u|^2\right)^{\frac{\beta}{2}}X^\e_\ell u\right)$ in the latter equation and, using the rule $$\sum_{i=1}^{2n+1}D_{i,j}^2f_\delta(\nabla_\e u)X^\e_iX^\e_\ell u=\sum_{i=1}^{2n+1}D_{i,j}^2f_\delta(\nabla_\e u)X^\e_\ell X^\e_i-s_\ell D^2_{\ell+s_\ell n,j}f_\delta(\nabla_\e u) Zu,$$
we obtain that for every $\ell=1,\ldots,2n$ it holds
\begin{equation}\label{4marzo25_1}
\aligned
    &\frac{1}{2} \int_0^T\int_\Om \eta^2\left(\delta+|\nabla_\e u|^2\right)^{\frac{\beta}{2}} \p_t\left[ (X^\e_\ell u)^2\right]\ dxdt\\&+\int_0^T\int_\Om\eta^2\left(\delta+|\nabla_\e u|^2\right)^{\frac{\beta}{2}}\sum_{i,j=1}^{2n+1}\Dfdeltaij(\nabla_\e u)X^\e_\ell X^\e_j  uX^\e_\ell X^\e_i u\ dxdt
    \\&+\int_0^T\int_\Om\eta^2 X^\e_\ell u \sum_{i,j=1}^{2n+1}\Dfdeltaij(\nabla_\e u)X^\e_\ell X^\e_j  u X^\e_i\left(\left(\delta+|\nabla_\e u|^2\right)^{\frac{\beta}{2}}\right)\ dxdt
    \\&=-2\int_0^T\int_\Om \eta X^\e_\ell u (\delta+|\nabla_\e u|^2)^{\frac{\beta}{2}}\sum_{i,j=1}^{2n+1}\Dfdeltaij (\nabla_\e u)X^\e_\ell X^\e_j u  X^\e_i\eta \ dxdt \\&+s_\ell\int_0^T\int_\Om\eta^2 \left(\delta+|\nabla_\e u|^2\right)^{\frac{\beta}{2}}\sum_{j=1}^{2n+1}D^2_{\ell+s_\ell n,j}f_\delta(\nabla_\e u)X^\e_\ell X^\e_ju Zu\ dxdt
    \\&-s_\ell\int_0^T\int_\Om \Dfdeltanl(\nabla_\e u) Z\left(\eta^2(\delta+|\nabla_\e u|^2)^{\frac{\beta}{2}}X^\e_\ell u\right)\ dxdt    =I^1_\ell+I^2_\ell+I^3_\ell.
\endaligned
\end{equation}

If $\ell=2n+1$, we use $\phi = \eta^2\left(\delta+|\nabla_\e u|^2\right)^{\frac{\beta}{2}}X^\e_\ell u$ as test function in the weak formulation of \eqref{eqX2n+1}. Again integrating by parts and computing the derivative we obtain
\begin{equation}\label{4marzo25_2}
\aligned
    &\frac{1}{2} \int_0^T\int_\Om \eta^2\left(\delta+|\nabla_\e u|^2\right)^{\frac{\beta}{2}} \p_t\left[ (X^\e_{2n+1} u)^2\right]\ dxdt\\&+\int_0^T\int_\Om\eta^2\left(\delta+|\nabla_\e u|^2\right)^{\frac{\beta}{2}}\sum_{i,j=1}^{2n+1}\Dfdeltaij(\nabla_\e u)X^\e_{2n+1} X^\e_j  uX^\e_{2n+1} X^\e_i u\ dxdt
    \\&+\int_0^T\int_\Om\eta^2 X^\e_{2n+1} u \sum_{i,j=1}^{2n+1}\Dfdeltaij(\nabla_\e u)X^\e_{2n+1} X^\e_j  u X^\e_i\left(\left(\delta+|\nabla_\e u|^2\right)^{\frac{\beta}{2}}\right)\ dxdt
    \\&=-2\int_0^T\int_\Om \eta X^\e_{2n+1} u (\delta+|\nabla_\e u|^2)^{\frac{\beta}{2}}\sum_{i,j=1}^{2n+1}\Dfdeltaij (\nabla_\e u)X^\e_{2n+1} X^\e_j u  X^\e_i\eta \ dxdt=:I^1_{2n+1}.
\endaligned
\end{equation}

Since the left hand side is the same for all values of $\ell=1,\ldots,2n+1$,
we will handle together the last integral in it: computing the derivative and using again the rule
\begin{equation}\label{18aprile25}
    \sum_{k=1}^{2n+1}X^\e_iX^\e_kuX^\e_ku=\sum_{k=1}^{2n+1}X^\e_k X^\e_iuX^\e_ku+s_i ZuX^\e_{i+s_in}u,
\end{equation}
we obtain
\begin{equation}\label{4marzo25_3}
\aligned
	&\int_0^T\int_\Om\eta^2 X^\e_\ell u \sum_{i,j=1}^{2n+1}\Dfdeltaij(\nabla_\e u)X^\e_\ell X^\e_j  u X^\e_i\left(\left(\delta+|\nabla_\e u|^2\right)^{\frac{\beta}{2}}\right)\ dxdt\\&=\beta\int_\Omega \eta^2\left(\delta+|\nabla_\e u|^2\right)^\frac{\beta-2}{2}\sum_{i,j=1}^{2n+1}{\Dfdeltaij} (\nabla_\e u)X^\e_\ell X^\e_juX^\e_\ell u\sum_{k=1}^{2n+1}X^\e_kX_i^\e uX^\e_kudx\\&+s_i\beta \int_\Omega \eta^2\left(\delta+|\nabla_\e u|^2\right)^\frac{\beta-2}{2}X^\e_\ell u Zu\sum_{i,j=1}^{2n} {\Dfdeltaij} (\nabla_\e u)X^\e_\ell X^\e_ju X^\e_{i+s_in}udx
\endaligned
\end{equation}
where we recall that $s_i=(-1)^{[\frac{i}{n+1}]}$, $i=1,\ldots,2n$.

We denote by $I^\ell_4$ the last integral in the right hand side of \eqref{4marzo25_3} and we set $I^2_{2n+1}=I^3_{2n+1}=0$. Using \eqref{4marzo25_3} in \eqref{4marzo25_1} and \eqref{4marzo25_2} and summing up over $\ell=1,\ldots,2n+1$ we obtain:
\begin{equation*}
    \aligned
    &\frac{1}{\beta+2} \int_0^T\int_\Om\eta^2 \p_t \left[(\delta+|\nabla_\e u|^2)^{\frac{\beta}{2}+1}\right]\ dxdt \\&+\int_0^T\int_\Om\eta^2 (\delta+|\nabla_\e u|^2)^{\frac{\beta}{2}}\sum_{i,j,\ell=1}^{2n+1} \Dfdeltaij (\nabla_\e u) X^\e_\ell X^\e_j  u X^\e_\ell X^\e_i u\ dxdt\\ &+\beta\int_\Omega \eta^2\left(\delta+|\nabla_\e u|^2\right)^\frac{\beta-2}{2}\sum_{i,j=1}^{2n+1}{\Dfdeltaij} (\nabla_\e u)\sum_{\ell=1}^{2n+1}X^\e_\ell X^\e_juX^\e_\ell u\sum_{k=1}^{2n+1}X^\e_kX_i^\e uX^\e_kudx\\&=\sum_{\ell=1}^{2n+1}\left(I_1^\ell+I_2^\ell+I_3^\ell-I_4^\ell\right)
    \endaligned
\end{equation*}
Using the structure condition \eqref{structurecondition01} and the fact that the third term in the left-hand side is always positive, it follows that  
\begin{equation}\label{4marzo25_4}
\aligned
    &\frac{1}{\beta+2} \int_0^T\int_\Om\eta^2 \p_t \left[(\delta+|\nabla_\e u|^2)^{\frac{\beta}{2}+1}\right]\ dxdt\\&+\int_0^T\int_\Om\eta^2 (\delta+|\nabla_\e u|^2)^{\frac{\beta}{2}}\sum_{i=1}^{n}\lambda_{i,\delta}(\nabla_\e u)\left(|\nabla_\e X^\e_i u|^2+|\nabla_\e X^\e_{i+n} u|^2+|\nabla_\e X^\e_{2n+1} u|^2\right)\ dxdt\\ &\leq \sum_{\ell=1}^{2n+1}\left(I_1^\ell+I_2^\ell+I_3^\ell+|I_4|^\ell\right).
\endaligned
\end{equation}
The structure condition \eqref{structurecondition01} and the Young inequality imply
\begin{equation}\label{gnam1}
\aligned
    &\sum_{\ell=1}^{2n+1}I^1_\ell\leq 2\int_0^T\int_\Om \eta (\delta+|\nabla_\e u|^2)^{\frac{\beta+1}{2}}\sum_{\ell=1}^{2n+1}\left|\sum_{i,j=1}^{2n+1}\Dfdeltaij(\nabla_\e u)X^\e_\ell X^\e_ju   X^\e_i\eta\right| \ dxdt\\& \leq \tau \int_0^T\int_\Om\eta^2 (\delta+|\nabla_\e u|^2)^{\frac{\beta}{2}}\sum_{i=1}^{n}\lambda_{i,\delta}(\nabla_\e u)\left(|\nabla_\e X^\e_iu|^2+|\nabla_\e X^\e_{i+n}u|^2+|\nabla_\e X^\e_{2n+1}u|^2\right)\ dxdt\\ &+ \frac{c}{\tau}\int_0^T\int_\Om |\nabla_\e \eta|^2(\delta+|\nabla_\e u|^2)^{\frac{\beta+2}{2}}\sum_{i=1}^n\lambda_{i,\delta}(\nabla_\e u) \ dxdt.
\endaligned
\end{equation}
As for $I_2^\ell$, again the structure condition \eqref{structurecondition01} and the Young inequality implies that
\begin{equation}\label{gnam2}
    \aligned
    &\sum_{\ell=1}^{2n+1}I^2_\ell
    \le \tau \int_0^T\int_\Om\eta^2 (\delta+|\nabla_\e u|^2)^{\frac{\beta}{2}}\sum_{i=1}^{n}\lambda_{i,\delta}(\nabla_\e u)\left(|\nabla_\e X^\e_i u|^2+|\nabla_\e X^\e_{i+n} u|^2+|\nabla_\e X^\e_{2n+1} u|^2\right)\ dxdt\\& +\frac{c}{\tau}\int_0^T\int_\Om \eta^2 (\delta+|\nabla_\e u|^2)^{\frac{\beta}{2}} \sum_{i=1}^n\lambda_{i,\delta}(\nabla_\e u)|Z u|^2\ dxdt.
    \endaligned
\end{equation}
Let us estimate $I^3_\ell$, for any $\ell\in \{1,\ldots,2n+1\}$. Computing the derivative and integrating by parts, one has
\begin{equation*}
\aligned
    I^3_\ell = &-2s_\ell \int_0^T\int_\Om  \eta Z\eta\Dfdeltanl(\nabla_\e u)(\delta+|\nabla_\e u|^2)^{\frac{\beta}{2}} X^\e_\ell u\ dxdt\\& -s_\ell\beta\int_0^T\int_\Om \eta^2 X^\e_\ell u\Dfdeltanl(\nabla_\e u)  (\delta+|\nabla_\e u|^2)^{\frac{\beta-2}{2}} \sum_{k=1}^{2n+1}X^\e_k u X^\e_k Zu \ dxdt\\& -s_\ell\int_0^T\int_\Om \eta^2 (\delta+|\nabla_\e u|^2)^{\frac{\beta}{2}}\Dfdeltanl(\nabla_\e u)X^\e_\ell Zu\ dxdt \\ = &-2s_\ell \int_0^T\int_\Om  \eta Z\eta\Dfdeltanl(\nabla_\e u)(\delta+|\nabla_\e u|^2)^{\frac{\beta}{2}} X^\e_\ell u\ dxdt\\&+s_\ell\beta\sum_{k=1}^{2n+1}\int_0^T\int_\Om X^\e_k\left(\eta^2 X^\e_\ell u\Dfdeltanl(\nabla_\e u)  (\delta+|\nabla_\e u|^2)^{\frac{\beta-2}{2}} X^\e_k u\right) Zu \ dxdt\\& +s_\ell\int_0^T\int_\Om X^\e_\ell\left(\eta^2 (\delta+|\nabla_\e u|^2)^{\frac{\beta}{2}}\Dfdeltanl(\nabla_\e u)\right) Zu\ dxdt.
\endaligned
\end{equation*} 

Therefore, computing derivatives and using commutation rules analogous to \eqref{18aprile25} in the last two integrals, using the structure condition \eqref{structurecondition01}, one has
\begin{equation}\label{gnam3}
\aligned
    \sum_{\ell=1}^{2n+1}I^3_\ell\leq &\tau \int_0^T\int_\Om\eta^2 (\delta+|\nabla_\e u|^2)^{\frac{\beta}{2}}\sum_{i=1}^{n}\lambda_{i,\delta}(\nabla_\e u)\left(|\nabla_\e X^\e_i u|^2+|\nabla_\e X^\e_{i+n} u|^2+|\nabla_\e X^\e_{2n+1} u|^2\right)\ dxdt\\ & + c(\beta + 1)\int_0^T\int_\Om \left(|\nabla_\e\eta|^2+\eta|Z\eta|\right)(\delta+|\nabla_\e u|^2)^{\frac{\beta+2}{2}}\sum_{i=1}^{n}\lambda_{i,\delta}(\nabla_\e u)\ dxdt \\& +\frac{c(\beta +1)^2}{\tau}\int_0^T\int_\Om \eta^2 (\delta+|\nabla_\e u|^2)^{\frac{\beta}{2}} \sum_{i=1}^{n}\lambda_{i,\delta}(\nabla_\e u)|Z u|^2\ dxdt.
\endaligned
\end{equation}

In the end,
\begin{equation}\label{gnam4}
\aligned
    \sum_{\ell=1}^{2n+1}|I^4_\ell|\leq &\tau\int_0^T\int_\Om\eta^2 (\delta+|\nabla_\e u|^2)^{\frac{\beta}{2}}\sum_{i=1}^{n}\lambda_{i,\delta}(\nabla_\e u)\left(|\nabla_\e X^\e_i u|^2+|\nabla_\e X^\e_{i+n} u|^2+|\nabla_\e X^\e_{2n+1} u|^2\right)\ dxdt\\ & +\frac{c(\beta+1)^2}{\tau}\int_0^T\int_\Om \eta^2 (\delta+|\nabla_\e u|^2)^{\frac{\beta}{2}} \sum_{i=1}^{n}\lambda_{i,\delta}(\nabla_\e u)|Z u|^2\ dxdt.
\endaligned   
\end{equation}

The thesis follows putting together \eqref{4marzo25_4}, \eqref{gnam1}, \eqref{gnam2}, \eqref{gnam3} and \eqref{gnam4}, choosing for instance $\tau=\frac{1}{8}$ and using the Leibniz rule on the first term in the left hand side, eventually redefining the constant $c=c(n,p,L)$.
\end{proof}

\section{An integrability estimate for the vertical derivative of the approximating solutions}\label{S:integrabilityZu}

In this section, we establish a local uniform integrability estimate for the vertical derivative of the approximating solutions $Z u$, stated in Proposition \ref{Poincaretypeineq}. This result will be used in the next section to remove the term involving $Z u$ in Lemma \ref{CaccioppoliXu}, see Theorem \ref{maincaccioppoli}.

The first step toward proving Lemma \ref{Poincaretypeineq} is to derive the following standard, uniform Caccioppoli-type inequality for $Z u$.

\begin{lemma} 
Let $u$ be a weak solution to \eqref{approx1} in $Q$.
For any $\beta \geq 0$ and any non-negative function $\eta\in C^{1}([0,T], C^\infty_c (\Om))$ vanishing on the parabolic boundary of $Q$, one has
\begin{equation}\label{caccioppoliZu}
\aligned
    &\int_0^T\int_\Omega\eta^2|Zu|^{\beta}\sum_{i=1}^n \lambda_{i,\delta}(\nabla_\e u)\left((X^\e_iZu)^2+(X^\e_{i+n}Zu)^2+(X^\e_{2n+1}Zu)^2\right)\ dxdt\\&\leq\frac{2L^2}{(\beta + 1)^2}\int_0^T\int_\Omega|\nabla_\e\eta|^2|Zu|^{\beta+2}\sum_{i=1}^n \lambda_{i,\delta}(\nabla_\e u)\ dxdt+\frac{2}{(\beta + 1)^2}\int_0^T\int_\Omega\eta|\partial_t\eta| |Zu|^{\beta+2}\ dxdt.
\endaligned
\end{equation}
\end{lemma}

\begin{proof}
We use $\phi= \eta^2 |Zu|^{\beta}Zu$ as a test function in the equation satisfied by $Zu$, i.e. \eqref{eqX2n+1} with $\e=1$, to obtain 
\begin{align}\label{14feb}
    \int_0^T\int_\Omega \partial_t Zu \eta^2 |Zu|^{\beta}Zu\ dxdt =\int_0^T\int_\Omega \sum_{i,j=1}^{2n+1}X^\e_i(\Dfdeltaij(\nabla_\e u)X^\e_j Zu) \eta^2 |Zu|^{\beta}Zu\ dxdt.
\end{align}
Integrating by parts, the left-hand side of \eqref{14feb} can be expressed as 
$$\int_0^T\int_\Omega \eta^2\partial_t Zu |Zu|^{\beta}Zu\ dxdt=
\frac{1}{\beta + 2}\int_0^T\int_\Omega \eta^2\partial_t \left(|Zu|^{\beta+2} \right)\ dxdt=-\frac{2}{\beta + 2}\int_0^T\int_\Omega\eta\p_t\eta |Zu|^{\beta+2}\ dxdt. 
$$
An integration by parts in the right-hand side of \eqref{14feb} gives
\begin{multline*}
    \int_0^T\int_\Omega \sum_{i,j=1}^{2n+1}X^\e_i(\Dfdeltaij(\nabla_\e u)X^\e_j Zu) \eta^2 |Zu|^{\beta}Zu\ dxdt=\\- \int_0^T\int_\Omega \sum_{i,j=1}^{2n+1}\Dfdeltaij(\nabla_\e u)X^\e_j Zu X^\e_i(\eta^2 |Zu|^{\beta}Zu)\ dxdt\\= - 2\int_0^T\int_\Omega \eta|Zu|^{\beta}Zu\sum_{i,j=1}^{2n+1}\Dfdeltaij(\nabla_\e u)X^\e_j Zu X^\e_i\eta \ dxdt\\ -(\beta +1)\int_0^T\int_\Omega\eta^2 |Zu|^{\beta} \sum_{i,j=1}^{2n+1}\Dfdeltaij(\nabla_\e u)X^\e_j Zu  X^\e_i Zu\ dxdt.
\end{multline*}

Combining the previous equations we obtain
\begin{multline*}
\aligned
   &(\beta+1)\int_0^T\int_\Omega\eta^2 |Zu|^{\beta} \sum_{i,j=1}^{2n+1}\Dfdeltaij(\nabla_\e u)X^\e_j Zu  X^\e_i Zu\ dxdt\\&=- 2\int_0^T\int_\Omega \eta|Zu|^{\beta}Zu\sum_{i,j=1}^{2n+1}\Dfdeltaij(\nabla_\e u)X^\e_j Zu X^\e_i\eta \ dxdt
    +\frac{2}{\beta + 2}\int_0^T\int_\Omega\eta\p_t\eta |Zu|^{\beta+2}\ dxdt.
\endaligned
\end{multline*}
The structure condition \eqref{structurecondition01} and the Young inequality imply
\begin{multline*}
\aligned
    &\int_0^T\int_\Omega\eta^2|Zu|^{\beta}\sum_{i=1}^n \lambda_{i,\delta}(\nabla_\e u)\left((X^\e_iZu)^2+(X^\e_{i+n}Zu)^2+(X^\e_{2n+1}Zu)^2\right)\ dxdt\\ &\leq\frac{2}{(\beta + 1)}\int_0^T\int_\Omega\eta |Zu|^{\beta +1}\left|\sum_{i,j=1}^{2n+1}\Dfdeltaij(\nabla_\e u)X^\e_j Zu X^\e_i\eta\right| dxdt+\frac{2} {(\beta + 1)^2}\int_0^T\int_\Omega\eta|\partial_t\eta| |Zu|^{\beta+2}\ dxdt\\
    &\leq\frac{\tau}{2}\int_0^T\int_\Omega\eta^2 |Zu|^{\beta}\sum_{i=1}^n \lambda_{i,\delta}(\nabla_\e u)\left((X^\e_iZu)^2+(X^\e_{i+n}Zu)^2+(X^\e_{2n+1}Zu)^2\right)\ dxdt\\&+\frac{2L^2}{\tau(\beta + 1)^2}\int_0^T\int_\Omega|\nabla_\e\eta|^2|Zu|^{\beta+2}\sum_{i=1}^n \lambda_{i,\delta}(\nabla_\e u)\ dxdt+\frac{2}{(\beta + 1)^2}\int_0^T\int_\Omega\eta|\partial_t\eta| |Zu|^{\beta+2}\ dxdt,  
\endaligned
\end{multline*}
where in the last inequality we used the structure condition \eqref{structurecondition01} and the Young inequality. The thesis follows by choosing $\tau=1$.
\end{proof}

As already anticipated in the Introduction, a key result for the proof of Proposition \ref{Poincaretypeineq} is the following Poincaré-type inequality, established in \cite[Lemma 4.1]{capogna2021lipschitz}. This result is independent of the specific equation \eqref{maineq01} and holds for any function $u \in C^2_{\text{loc}}(Q)$. It is also the only point in the paper where the restriction $2 \leq p \leq 4$ is required.

Although we do not use this result directly, a crucial estimate from its proof plays an essential role in our argument. However, we state the inequality (in a slightly different form) for the reader's convenience.

\begin{lemma}\label{poincarétypeineq0}
Let $2\leq p\leq 4$ and $u\in C^2_{loc}(Q)$, and let us denote by $\nabla_{H,\ell}=(X^\e_\ell, X^\e_{n+\ell})$, for any $\ell\in\{1,\ldots, 2n\}$. There is a constant $c=c(n,p)$ such that, for any $\beta\geq0$ and any non negative $\eta\in C^1([0,T], C_c^\infty(\Omega))$ vanishing on the parabolic boundary of $Q$, we have
\begin{equation}\label{poincarétypeineq1}
    \aligned
    \int_0^T\int_\Om \eta^{\beta+p}|Zu|^{\beta+p}\ dxdt&\leq c(\beta+p)\|\nabla_H\eta\|_\infty\int\int_{\textnormal{supp}(\eta)} (\delta+ |\nabla_{H,\ell} u|^2)^{\frac{\beta+p}{2}}\ dxdt\\
    &+c(\beta+p)\int_0^T\int_\Omega\eta^{\beta+4}|Zu|^{\beta}(\delta + |\nabla_{H,\ell} u|^2)^{\frac{p-2}{2}}|\nabla_{H,\ell}Zu|^2\ dxdt.
    \endaligned
\end{equation}
\end{lemma}

\begin{prop}\label{Poincaretypeineq}
Let $2\leq p\leq 4$ and let $u$ be a weak solution to \eqref{approx1}. Then, 
\begin{equation*}
    Zu\in L^p_{loc}(Q).
\end{equation*}
Moreover, there is a constant $c=c(n,p,L)$ such that, for any $\beta\geq0$ and any non negative $\eta\in C^1([0,T], C_c^\infty(\Omega))$ vanishing on the parabolic boundary of $Q$, we have
\begin{multline}\label{mainestimateZu}
    \left(\int_0^T\int_\Om \eta^{\beta+p}|Zu|^{\beta+p}\ dxdt\right)^{\frac{1}{\beta+p}}\leq c(\beta+p) \|\nabla_\e\eta\|_\infty\left(\int\int_{\textnormal{supp}(\eta)} (\delta+ |\nabla_\e u|^2)^{\frac{\beta+p}{2}}\ dxdt\right)^{\frac{1}{\beta+p}}\\+c(\beta+p)\|\eta\p_t\eta\|_\infty^{\frac{1}{2}}\ |\textnormal{supp}(\eta)|^{\frac{p-2}{2(\beta+p)}}\left(\int\int_{\textnormal{supp}(\eta)} \left(\delta+|\nabla_\e u|^2\right)^{\frac{\beta+p}{2}}\ dxdt\right)^{\frac{4-p}{2(\beta+p)}}. 
\end{multline}
\end{prop}

\begin{proof}
If we denote by $\nabla_{H,\ell}=(X^\e_\ell, X^\e_{n+\ell})$, for any $\ell\in\{1,\ldots, 2n\}$, and by
\begin{equation*}\label{definteg}
    \aligned
    I&:=\int_0^T\int_\Om \eta^{\beta+p}|Zu|^{\beta+p}\ dxdt\\
    R&:=\int\int_{\textnormal{supp}(\eta)} (\delta+ |\nabla_{H,\ell} u|^2)^{\frac{\beta+p}{2}}\ dxdt\\
    M&:=\int_0^T\int_\Omega\eta^{4+\beta}|Zu|^{\beta}(\delta + |\nabla_{H,\ell} u|^2)^{\frac{p-2}{2}}|\nabla_{H,\ell}Zu|^2\ dxdt,
    \endaligned 
\end{equation*}
then from the proof of Lemma \ref{poincarétypeineq0} it follows that 
\begin{equation}\label{estimateZu1}
    I\leq 2(\beta+p)\left(M^{\frac{1}{2}}R^{\frac{4-p}{2(\beta+p)}}I^{\frac{2p-4+\beta}{2(\beta+p)}}+\|\nabla_{H,\ell}\eta\|_{L^\infty}R^{\frac{1}{\beta+p}}I^{\frac{\beta+p-1}{\beta+p}}\right).
\end{equation}
First of all we observe that 
\begin{equation}\label{estimateR}
       R\leq\int\int_{\textnormal{supp}(\eta)} (\delta+ |\nabla_\e u|^2)^{\frac{\beta+p}{2}}\ dxdt.
\end{equation}
On the other hand, we apply the inequality \eqref{caccioppoliZu} to estimate the integral $M$ in the following way:
\begin{equation}\label{estimateZu2}
\aligned
    M&\leq \int_0^T\int_\Omega\left(\eta^\frac{4+\beta}{2}\right)^2|Zu|^{\beta}\sum_{i=1}^n \lambda_{i,\delta}(\nabla_\e u)\left((X^\e_iZu)^2+(X^\e_{i+n}Zu)^2+(X^\e_{2n+1}Zu)^2\right)\ dxdt \\
    &\leq c\int_0^T\int_\Omega\eta^{2+\beta}|\nabla_\e\eta|^2|Zu|^{\beta+2}\sum_{i=1}^n \lambda_{i,\delta}(\nabla_\e u)\ dxdt  \\&+c\int_0^T\int_\Omega\eta^{3+\beta}|\partial_t\eta| |Zu|^{\beta+2}\ dxdt \\
    &\leq c\|\nabla_\e\eta\|^2_\infty\ I^\frac{\beta+2}{\beta+p} \left(\int\int_{\textnormal{supp}(\eta)} \left(\sum_{i=1}^n \lambda_{i,\delta}(\nabla_\e u)\right)^{\frac{\beta+p}{p-2}}\ dxdt\right)^{\frac{p-2}{\beta+p}}\\&+c\|\eta\p_t\eta\|_\infty\ I^{\frac{\beta+2}{\beta+p}} |\text{supp}(\eta)|^{\frac{p-2}{\beta+p}} \\
    &\leq c\|\nabla_\e\eta\|^2_\infty\ I^\frac{\beta+2}{\beta+p} \left(\int\int_{\textnormal{supp}(\eta)} \left(\delta+|\nabla_\e u|^2\right)^{\frac{\beta+p}{2}}\ dxdt\right)^{\frac{p-2}{\beta+p}}\\&+c\|\eta\p_t\eta\|_\infty\ |\text{supp}(\eta)|^{\frac{p-2}{\beta+p}}\ I^{\frac{\beta+2}{\beta+p}},
\endaligned
\end{equation}
where we used the H\"older inequality and the fact that
\begin{equation}\label{6marzo25_4}
    \sum_{i=1}^n\lambda_{i,\delta}(\nabla_\e u)\leq n\left(\delta+|\nabla_\e u|^2\right)^\frac{p-2}{2},
\end{equation}
eventually changing the constant $c=c(n,p,L)$. From \eqref{estimateZu1}, using the estimates \eqref{estimateR} and \eqref{estimateZu2}, together with the fact that $\|\nabla_{H,\ell}\eta\|_\infty\leq \|\nabla_\e\eta\|_\infty$, we can conclude:
\begin{multline*}
    I\leq c(\beta+p)\Bigg{(}\|\nabla_\e\eta\|_\infty\ I^\frac{\beta+p-1}{\beta+p} \left(\int_0^T\int_{\text{supp}(\eta)} \left(\delta+|\nabla_\e u|^2\right)^{\frac{\beta+p}{2}}\ dxdt\right)^{\frac{1}{\beta+p}}\\+\|\eta\p_t\eta\|_\infty^{\frac{1}{2}}\ |\text{supp}(\eta)|^{\frac{p-2}{2(\beta+p)}}\ I^{\frac{\beta+p-1}{\beta+p}} \left(\int_0^T\int_{\text{supp}(\eta)} \left(\delta+|\nabla_\e u|^2\right)^{\frac{\beta+p}{2}}\ dxdt\right)^{\frac{4-p}{2(\beta+p)}}\Bigg{)}.
\end{multline*}
The thesis follows by dividing both sides by the quantity $I^{\frac{\beta+p-1}{\beta+p}}$.
\end{proof}

\section{Main Caccioppoli-type inequality for the first derivatives of approximating solutions}\label{S:maincaccioppoli}

The next theorem presents a standard Caccioppoli-type inequality for the Riemannian gradient of the approximating solutions $\nabla_\varepsilon u$, in which the vertical derivative $Z u$ no longer appears.

\begin{thrm}\label{maincaccioppoli}
Let $2\leq p\leq 4$ and let $u$ be a weak solution to \eqref{approx1}. There is a constant $c=c(n,p,L)$ such that, for any $\beta\geq0$ and any non-negative $\eta\in C^1([0,T], C_c^\infty(\Omega))$ vanishing on the parabolic boundary of $Q$, we have
\begin{equation}\label{standardcaccioppoli}
\aligned
    &\sup_{0<t<T}\int_\Om\eta^2\left(\delta+|\nabla_\e u|^2\right)^{\frac{\beta+2}{2}}\ dx\\&+\int_0^T\int_\Om\eta^2 (\delta+|\nabla_\e u|^2)^{\frac{\beta}{2}}\sum_{i=1}^{n}\lambda_{i,\delta}(\nabla_\e u)\left(|\nabla_\e X^\e_i u|^2+|\nabla_\e X^\e_{i+n} u|^2+|\nabla_\e X^\e_{2n+1} u|^2\right)\ dxdt\\&\leq c(\beta+p)^5\left(\|\nabla_\e\eta\|_\infty^2+\|\eta Z\eta\|_\infty\right)\int\int_{\textnormal{supp}(\eta)}\left(\delta+|\nabla_\e u|^2\right)^{\frac{\beta+p}{2}}\ dx\\&+c(\beta+p)^5\|\eta\p_t\eta\|_\infty |\textnormal{supp}(\eta)|^{\frac{p-2}{\beta+p}}\left(\int\int_{\textnormal{supp}(\eta)}\left(\delta+|\nabla_\e u|^2\right)^{\frac{\beta+p}{2}}\ dx\right)^{\frac{\beta+2}{\beta+p}}.
\endaligned
\end{equation}
\end{thrm}

\begin{proof}
Lemma \ref{CaccioppoliXu} implies that
\begin{equation*}
    \aligned
    &\sup_{0<t<T}\int_\Om\eta^2\left(\delta+|\nabla_\e u|^2\right)^{\frac{\beta+2}{2}}\ dx\\&+\int_0^T\int_\Om\eta^2 (\delta+|\nabla_\e u|^2)^{\frac{\beta}{2}}\sum_{i=1}^{n}\lambda_{i,\delta}(\nabla_\e u)\left(|\nabla_\e X^\e_i u|^2+|\nabla_\e X^\e_{i+n} u|^2+|\nabla_\e X^\e_{2n+1} u|^2\right)\ dxdt\\
    &\leq c(\beta+p)^2\int_0^T\int_\Om \left(|\nabla_\e\eta|^2+\eta|Z\eta|\right)(\delta+|\nabla_\e u|^2)^{\frac{\beta+2}{2}}\sum_{i=1}^n\lambda_{i,\delta}(\nabla_\e u)\ dxdt\\
    &+c\int_0^T\int_\Om\eta|\p_t\eta|(\delta+|\nabla_\e u|^2)^{\frac{\beta+2}{2}}\ dx\\&+c(\beta+p)^3\int_0^T\int_\Om \eta^2 (\delta+|\nabla_\e u|^2)^{\frac{\beta}{2}}\sum_{i=1}^n\lambda_{i,\delta}(\nabla_\e u)|Z u|^2\ dxdt=:I_1+I_2+I_3.
    \endaligned
\end{equation*}
It is obvious that $I_1$ is bounded by the first term in right hand side of \eqref{standardcaccioppoli}. Hence, the thesis follows if the other two terms are also bounded by the right hand side of \eqref{standardcaccioppoli}. Let us estimate $I_2$: the H\"older's inequality gives us
\begin{equation*}
    \int_0^T\int_\Om\eta|\p_t\eta|(\delta+|\nabla_\e u|^2)^{\frac{\beta+2}{2}}\ dx\leq \|\eta\p_t\eta\|_\infty |\textnormal{supp}(\eta)|^{\frac{p-2}{\beta+p}}\left(\int\int_{\textnormal{supp}(\eta)}\left(\delta+|\nabla_\e u|^2\right)^{\frac{\beta+p}{2}}\ dx\right)^{\frac{\beta+2}{\beta+p}},
\end{equation*}
which is obviously bounded by the second term in the right hand side of \eqref{standardcaccioppoli}. In the end, for $I_3$ we use the H\"older inequality, \eqref{mainestimateZu} and \eqref{6marzo25_4} to obtain
\begin{equation*}
\aligned
    &\int_0^T\int_\Om \eta^2 (\delta+|\nabla_\e u|^2)^{\frac{\beta}{2}}\sum_{i=1}^n\lambda_{i,\delta}(\nabla_\e u)|Z u|^2\ dxdt\\&\leq \left(\int_0^T\int_\Om \eta^{\beta+p}|Z u|^{\beta+p}\ dxdt\right)^{\frac{2}{\beta+p}}\left(\int\int_{\textnormal{supp}(\eta)}\left((\delta+|\nabla_\e u|^2)^{\frac{\beta}{2}}\sum_{i=1}^n\lambda_{i,\delta}(\nabla_\e u)\right)^{\frac{\beta+p}{\beta+p-2}} \ dxdt\right)^{\frac{\beta+p-2}{\beta+p}}\\
    &\leq c\left(\int_0^T\int_\Om \eta^{\beta+p}|Z u|^{\beta+p}\ dxdt\right)^{\frac{2}{\beta+p}}\left(\int\int_{\textnormal{supp}(\eta)} (\delta+|\nabla_\e u|^2)^{\frac{\beta+p}{2}} \ dxdt\right)^{\frac{\beta+p-2}{\beta+p}}\\
    &\leq c(\beta+p)^2\|\nabla_\e \eta\|_\infty^2\int\int_{\textnormal{supp}(\eta)}\left(\delta+|\nabla_\e u|^2\right)^{\frac{\beta+p}{2}}\ dxdt\\&+    c(\beta+p)^2\|\eta\p_t\eta\|_\infty|\textnormal{supp}(\eta)|^{\frac{p-2}{\beta+p}}\left(\int\int_{\textnormal{supp}(\eta)}\left(\delta+|\nabla_\e u|^2\right)^{\frac{\beta+p}{2}}\ dxdt\right)^{\frac{\beta+2}{\beta+p}},
\endaligned
\end{equation*}
and the thesis follows.
\end{proof}

The following result is a straightforward consequence of the previous theorem. The inequality \eqref{standardcaccioppoli2} will serve as a key ingredient, together with the Sobolev inequality stated in Lemma \ref{lemma-sobolev}, to implement the Moser iteration scheme.

\begin{cor}
Let $2\leq p\leq 4$ and let $u$ be a weak solution to \eqref{approx1}. There is a constant $c=c(n,p,L)$ such that, for any $\beta\geq0$ and any non-negative $\eta\in C^1([0,T], C_c^\infty(\Omega))$ vanishing on the parabolic boundary of $Q$, we have\begin{equation}\label{standardcaccioppoli2}
\aligned
    &\sup_{0<t<T}\int_\Om\eta^2\left(\delta+(X^\e_\ell u)^2\right)^{\frac{\beta+2}{2}}\ dx+\int_0^T\int_\Om\eta^2 (\delta+(X^\e_\ell u)^2)^{\frac{\beta+p-2}{2}}|\nabla_\e X^\e_\ell u|^2\ dxdt\\&\leq c(\beta+p)^5\left(\|\nabla_\e\eta\|_\infty^2+\|\eta Z\eta\|_\infty\right)\int\int_{\textnormal{supp}(\eta)}\left(\delta+|\nabla_\e u|^2\right)^{\frac{\beta+p}{2}}\ dx\\&+c(\beta+p)^5\|\eta\p_t\eta\|_\infty |\textnormal{supp}(\eta)|^{\frac{p-2}{\beta+p}}\left(\int\int_{\textnormal{supp}(\eta)}\left(\delta+|\nabla_\e u|^2\right)^{\frac{\beta+p}{2}}\ dx\right)^{\frac{\beta+2}{\beta+p}},
\endaligned
\end{equation}
for every $\ell=1,\ldots,2n+1.$
\end{cor}

\section{Uniform Lipschitz estimate for approximating solutions}\label{S:lipschitzestimate}

In Theorem \ref{moser}, we establish a uniform local $L^\infty$ bound for the Riemannian gradient of the approximating solutions. The argument is based on the Moser iteration scheme and relies on the observation that the quantity $\delta + |\nabla_\varepsilon u|$ is bounded below by $\delta$, and that for every $\beta \geq 0$, it is locally bounded in $L^{\beta + p}$ over a parabolic cylinder, uniformly in $\varepsilon$ and $\delta$.

In the iteration process, we consider the Riemannian balls $B_\varepsilon$ introduced in Subsection \ref{SS:riemannianapprox}. We recall that these balls converge to the sub-Riemannian balls in the sense of Hausdorff distance, ensuring that the resulting estimates remain stable as $\varepsilon, \delta \to 0^+$.

We begin by recalling the definition of a parabolic cylinder.

\begin{dfn}\label{riemanniancylinder}
Let $(x_0,t_0)\in \Om\times (0,T)$. For $r,\mu>0$, the Riemannian parabolic cylinder $Q^\e_{\mu,r}(x_0, t_0)\subset Q$ of center $(x_0,t_0)$ is the set 
$$Q^\e_{\mu,r}(x_0, t_0)= B_{\e}(x_0, r) \times (t_0-\mu\,r, t_0 ).$$

We call parabolic boundary of the cylinder $Q^\e_{\mu,r}(x_0, t_0)\subset Q$
the set 
$$\p Q^\e_{\mu,r}(x_0, t_0):=B_{\e}(x_0, r) \times \{t_0-\mu\,r\} \cup \partial B_{\e}(x_0, r) \times [t_0-\mu\,r, t_0).$$
\end{dfn}

We also state the following inequality, which plays a crucial role in the Moser iteration. It is a direct consequence of the Riemannian Sobolev inequality.

\begin{lemma}\label{lemma-sobolev}
Let $v\in C^\infty(Q)$ such that, for any $0<t<T$, the function $v(\cdot, t)$ has compact support in $\Om\times\{t\}$. There is $c=c(n,p)>0$ such that, for any $\e \in (0,1]$, one has
\begin{equation}\label{lemma-sobolev1}
    \|v\|_{L^{\frac{2N}{N-2}, 2}(Q)} \leq c \|\nabla_\e v\|_{L^{2,2}(Q,\mathbb{R}^{2n+1})}.
\end{equation}
\end{lemma}

It is worth noting that, as $\varepsilon$ tends to zero, the underlying geometry transitions from Riemannian to sub-Riemannian. Nevertheless, the constant $c$ in Lemma \ref{lemma-sobolev} remains uniform with respect to $\varepsilon$, as shown in \cite{capogna2016regularity}.

\begin{thrm}\label{moser}
Let $2\leq p\leq4$ and $u$ be a weak solution to \eqref{approx1}, then
\begin{equation*}
    \nabla_\e u\in L_{loc}^\infty(Q,\mathbb{R}^{2n+1}).
\end{equation*}
Moreover, there is a constant $c=c(n,p,L)$ such that for any $Q^\e_{\mu,2r}\subset Q$ it holds 
\begin{equation}\label{moser-conclusion}
    \|\nabla_\e u\|_{L^\infty(Q^\e_{\mu,r},\mathbb{R}^{2n+1})}\leq c\mu^{\frac{1}{2}}\,\max\left\{\left(\frac{1}{\mu r^{N+2}}\int\int_{Q^\e_{\mu,2r}}\left(\delta+|\nabla_\e u|^2\right)^{\frac{p}{2}}\ dxdt\right)^{\frac{1}{2}},\mu^{\frac{p}{2(2-p)}}\right\},
\end{equation}
where $c=c(n,p,L)>0$.
\end{thrm}

\begin{proof} 
We consider a non negative cut-off function $\eta\in C^1\left([0,T], C^\infty_c(\Omega)\right)$,  vanishing on the parabolic boundary of $Q$, such that $|\eta|\leq1$ in $Q$. For any $\beta\geq0$, we denote by
\begin{equation*}
\aligned
    &v_\ell=\eta(\delta+(X^\e_\ell u)^2)^{\frac{\beta+p}{4}},\quad \ell=1,\ldots,2n+1,\\
    &v=\eta(\delta+|\nabla_\e u|^2)^{\frac{\beta+p}{4}}.
\endaligned
\end{equation*}
The Caccioppoli inequality \eqref{standardcaccioppoli2} implies
\begin{equation}\label{5marzo25_1}
\aligned
    \sup_{0<t<T}\int_\Om v_\ell^m\ dx+\int_0^T\int_\Om|\nabla_\e v_\ell|^2\ dxdt&\leq c(\beta+p)^5\left(\|\nabla_\e\eta\|_\infty^2+\|\eta Z\eta\|_\infty\right)\int\int_{\textnormal{supp}(\eta)}v^2\ dx\\&+c(\beta+p)^5\|\eta\p_t\eta\|_\infty |\textnormal{supp}(\eta)|^{\frac{p-2}{\beta+p}}\left(\int\int_{\textnormal{supp}(\eta)}v^2\ dx\right)^{\frac{\beta+2}{\beta+p}},  
\endaligned
\end{equation}
where $m=\frac{2(\beta+2)}{\beta+p}$. We note that $4/p<m\leq 2$.

Moreover, 
\begin{equation*}
\aligned
    \int_0^T\int_\Om v_\ell^q\ dxdt&\leq \int_0^T\left(\int_\Om v_\ell^m\ dx\right)^{\frac{2}{N}}\left(\int_\Om v_\ell^{\frac{2N}{N-2}}\ dx\right)^{\frac{N-2}{N}}\ dt\\&\leq c\left(\sup_{0<t<T}\int_\Om v_\ell^m\ dx\right)^{\frac{2}{N}}\left(\int_0^T\int_\Om |\nabla_\e v_\ell|^2\ dxdt\right),
\endaligned
\end{equation*}
where we denoted by $$q=\frac{2(m+N)}{N}=2+\frac{4(\beta+2)}{N(\beta+p)}$$ and in the second inequality we used the Sobolev inequality \eqref{lemma-sobolev1}, with $c=c(n)$. Raising both sides of the previous inequality to the power $\frac{N}{N+2}$, using the Young inequality, \eqref{5marzo25_1} and summing over $\ell=1,\ldots, 2n+1$, we get
\begin{equation}\label{5marzo25_2}
\aligned
    \sum_{\ell=1}^{2n+1}\left(\int_0^T\int_\Om v_\ell^q\ dxdt\right)^{\frac{N}{N+2}}&\leq c(\beta+p)^5\left(\|\nabla_\e\eta\|_\infty^2+\|\eta Z\eta\|_\infty\right)\int\int_{\textnormal{supp}(\eta)}v^2\ dxdt\\&+c(\beta+p)^5\|\eta\p_t\eta\|_\infty |\textnormal{supp}(\eta)|^{\frac{p-2}{\beta+p}}\left(\int\int_{\textnormal{supp}(\eta)}v^2\ dx\right)^{\frac{\beta+2}{\beta+p}}.
\endaligned
\end{equation}
On the other hand,
\begin{equation}\label{5marzo25_3}
    \left(\int_0^T\int_\Om v^q\ dxdt\right)^{\frac{N}{N+2}}\leq (2n+1)^{\beta+p}\sum_{\ell=1}^{2n+1}\left(\int_0^T\int_\Om v_\ell^q\ dxdt\right)^{\frac{N}{N+2}},
\end{equation}
then putting together \eqref{5marzo25_2} and \eqref{5marzo25_3} we get the following inequality that we will iterate:
\begin{equation}\label{5marzo25_4}
\aligned
    \left(\int_0^T\int_\Om v^q\ dxdt\right)^{\frac{N}{N+2}}&\leq c(\beta+p)^5\left(\|\nabla_\e\eta\|_\infty^2+\|\eta Z\eta\|_\infty\right)\int\int_{\textnormal{supp}(\eta)}v^2\ dxdt\\&+c(\beta+p)^5\|\eta\p_t\eta\|_\infty |\textnormal{supp}(\eta)|^{\frac{p-2}{\beta+p}}\left(\int\int_{\textnormal{supp}(\eta)}v^2\ dx\right)^{\frac{\beta+2}{\beta+p}}.
\endaligned
\end{equation}
We consider $Q^\e_{\mu,2r}\subset Q$ and we define a sequence of radii $r_i=(1+2^{-i})r$ and a sequence of exponents $\beta_i$, such that $\beta_0=0$ and
\begin{equation*}
    \beta_{i+1}+p=(\beta_i+p)\left(1+\frac{2(\beta_i+2)}{N(\beta_i+p)}\right),
\end{equation*}
that is,
\begin{equation*}
    \beta_i=2(k^i-1),\quad \textnormal{with }k=\frac{N+2}{N}.
\end{equation*}
We denote by $Q_i=Q^\e_{\mu,r_i}$, so that $Q_0=Q^\e_{\mu,2r}$ and $Q_\infty=Q^\e_{\mu,r}$. Moreover we choose a standard parabolic cut-off function $\eta_i\in C^\infty(Q_i)$ such that $\eta_i=1$ in $Q_{i+1}$ and 
\begin{equation*}
    |\nabla_\e\eta_i|\leq \frac{2^{i+8}}{r},\quad |Z\eta_i|\leq\frac{2^{2i+8}}{r^2},\quad |\p_t\eta_i|\leq\frac{2^{2i+8}}{\mu r^2} \quad\textnormal{in }Q_i.
\end{equation*}
Now we take $\eta=\eta_i$ and $\beta=\beta_i$ in \eqref{5marzo25_4} and we obtain
\begin{equation}\label{6marzo25_1}
\aligned
    \left(\int\int_{Q_{i+1}}\left(\delta+|\nabla_\e u|^2\right)^{\frac{\alpha_i+1}{2}}\ dxdt\right)^{\frac{N}{N+2}}&\leq c 2^{2i}\alpha_i^7r^{-2}\Bigg{[}\left(\int\int_{Q_i}\left(\delta+|\nabla_\e u|^2\right)^{\frac{\alpha_i}{2}}\ dxdt\right)^{\frac{p-2}{\alpha_i}}\\&+\mu^{-1}\left(\mu r^{N+2}\right)^{\frac{p-2}{\alpha_i}}\Bigg{]}\left(\int\int_{Q_i}\left(\delta+|\nabla_\e u|^2\right)^{\frac{\alpha_i}{2}}\ dxdt\right)^{\frac{\alpha_i-p+2}{\alpha_i}},
\endaligned
\end{equation}
where $c=c(n,p,L)>0$ and $\alpha_i=\beta_i+p=p-2+2k^i$. We denote by
\begin{equation*}
    M_i=\left(\frac{1}{\mu r^{N+2}}\int\int_{Q_i}\left(\delta+|\nabla_\e u|^2\right)^\frac{\alpha_i}{2}\ dxdt\right)^{\frac{1}{\alpha_i}}.
\end{equation*}
Then we can rewrite \eqref{6marzo25_1} as
\begin{equation*}
    M_{i+1}^{\frac{\alpha_i+1}{k}}\leq c\mu^{\frac{2}{N+2}}2^{2i}\alpha_i^7\left(M_i^{p-2}+\mu^{-1}\right)M_i^{\alpha_i-p+2}.
\end{equation*}
We set 
\begin{equation*}
    \overline{M}_i=\max\left\{M_i,\mu^{\frac{1}{2-p}}\right\},
\end{equation*}
so that the above inequality implies
\begin{equation}\label{6marzo25_2}
    \overline{M}_{i+1}^{\frac{\alpha_i+1}{k}}\leq c\mu^{\frac{2}{N+2}}2^{2i}\alpha_i^7\overline{M}_i^{\alpha_i}
\end{equation}
Iterating \eqref{6marzo25_2} we obtain
\begin{equation*}
    \overline{M}_{i+1}\leq \left(\Pi_{j=0}^iK_j^{\frac{k^{i+1-j}}{\alpha_{i+1}}}\right)\overline{M}_0^{\frac{\alpha_0k^{i+1}}{\alpha_{i+1}}},
\end{equation*}
where
\begin{equation*}
    K_i=c\mu^{\frac{2}{N+2}}2^{2i}\alpha_i^7.
\end{equation*}
Letting $i$ tend to $\infty$, we obtain
\begin{equation*}
    \sup_{Q^\e_{\mu,r}}|\nabla_\e u|\leq\overline{M}_\infty=\limsup_{i\to\infty}\overline{M}_i\leq c\mu^{\frac{1}{2}}\overline{M}_0^{\frac{p}{2}},
\end{equation*}
where
\begin{equation*}
    \overline{M}_0=\max\left\{\left(\frac{1}{\mu r^{N+2}}\int\int_{Q^\e_{\mu,2r}}\left(\delta+|\nabla_\e u|^2\right)^{\frac{p}{2}}\ dxdt\right)^{\frac{1}{p}},\mu^{\frac{1}{2-p}}\right\}
\end{equation*}
which concludes the proof.
\end{proof}

\section{Higher integrability for the time derivative of approximating solutions}\label{S:integtimederivative}

This section is devoted to establishing a uniform local bound for the $L^q$ norm of the time derivative of the approximating solutions $\partial_t u$, for any $q \geq 1$, see Proposition \ref{qintegrabilitytime}. 

To achieve this, we begin with a Caccioppoli-type inequality for $\partial_t u$, which will serve as the key tool in deriving the desired integrability estimate.

\begin{lemma}
Let $u$ be a weak solution of \eqref{approx1}. There is $c=c(n,p,L) > 0$ such that, for any $\beta \geq 0$ and any non-negative function $\eta\in C^{1}([0,T], C^\infty_c (\Om))$ vanishing on the parabolic boundary of $Q$, one has
\begin{equation}\label{caccioppolidtu}
\aligned
    &\int_0^T\int_\Om\eta^{\beta+4}|\p_tu|^\beta\sum_{i=1}^{n}\lambda_{i,\delta}(\nabla_\e u)\left((X^\e_i\p_tu)^2+(X^\e_{i+n}\p_tu)^2+(X^\e_{2n+1}\p_tu)^2\right)\ dxdt\\
    &\leq\int_0^T\int_\Om\eta^{\beta+2}|\p_tu|^\beta\sum_{i=1}^{n}\lambda_{i,\delta}(\nabla_\e u)|\nabla_\e\eta|^2\ dxdt\\&+c\int_0^T\int_\Om\eta^{\beta+3}|\p_t\eta|\ |\p_tu|^{\beta+2}\ dxdt.
\endaligned
\end{equation}
\end{lemma}

\begin{proof}
First of all we notice that the function $\p_tu$ solves the following equation:
\begin{equation}\label{eqpt}
    \p_t\left(\p_t u\right)=\sum_{i,j=1}^{2n+1}X^\e_i\left(\Dfdeltaij(\nabla_\e u)X^\e_j\p_t u\right).
\end{equation}
Let $\eta\in C^1\left([0,T],C_c^\infty(\Omega)\right)$ be a non-negative cut off function vanishing on the parabolic boundary of $Q$, we use $\phi=\eta^{\beta+4}|\p_tu|^\beta\p_tu$ as test function in the weak formulation of \eqref{eqpt}: integrating by parts on the right hand side and dividing both sides by $\beta+1$ we obtain
\begin{equation}\label{7marzo25_1}
\aligned
    &\int_0^T\int_\Om\eta^{\beta+4}|\p_tu|^\beta\sum_{i,j=1}^{2n+1}\Dfdeltaij(\nabla_\e u)X^\e_j\p_t uX^\e_i\p_tu\ dxdt\\&=-\frac{\beta+4}{\beta+1}\int_0^T\int_\Om\eta^{\beta+3}|\p_tu|^\beta\p_tu\sum_{i,j=1}^{2n+1}\Dfdeltaij(\nabla_\e u)X^\e_j\p_t uX^\e_i\eta\ dxdt\\&-\frac{1}{\beta+1}\int_0^T\int_\Om\eta^{\beta+4}|\p_tu|^\beta\p_tu\ \p_t\left( \p_tu\right)\ dxdt=I_1+I_2.
\endaligned
\end{equation}
Using the structure condition \eqref{structurecondition02}, the left hand side of \eqref{7marzo25_1} can be bounded from below by
\begin{equation*}
    (LHS)\geq\int_0^T\int_\Om\eta^{\beta+4}|\p_tu|^\beta\sum_{i=1}^{n}\lambda_{i,\delta}(\nabla_\e u)\left((X^\e_i\p_tu)^2+(X^\e_{i+n}\p_tu)^2+(X^\e_{2n+1}\p_tu)^2\right)\ dxdt.
\end{equation*}
Let us bound by above the right hand side of \eqref{7marzo25_1}. 

For the first term we use \eqref{structurecondition01} and the Young inequality:
\begin{equation*}
\aligned
    I_1&\leq \tau\int_0^T\int_\Om\eta^{\beta+4}|\p_tu|^\beta\sum_{i=1}^{n}\lambda_{i,\delta}(\nabla_\e u)\left((X^\e_i\p_tu)^2+(X^\e_{i+n}\p_tu)^2+(X^\e_{2n+1}\p_tu)^2\right)\ dxdt\\
    +&\frac{c}{\tau}\int_0^T\int_\Om\eta^{\beta+2}|\p_tu|^\beta\sum_{i=1}^{n}\lambda_{i,\delta}(\nabla_\e u)|\nabla_\e\eta|^2\ dxdt,
\endaligned
\end{equation*}
where $c=c(n,p,L)>0$.

Instead the second term in the right hand side can be handled in the following way
\begin{equation*}
\aligned
    I_2&=-\frac{1}{(\beta+1)(\beta+2)}\int_0^T\int_\Om\eta^{\beta+4}\ \p_t\left(|\p_tu|^{\beta+2}\right)\ dxdt\\&=\frac{\beta+4}{(\beta+1)(\beta+2)}\int_0^T\int_\Om\eta^{\beta+3}\p_t\eta\ |\p_tu|^{\beta+2}\ dxdt\\
    &\leq c_0\int_0^T\int_\Om\eta^{\beta+3}|\p_t\eta|\ |\p_tu|^{\beta+2}\ dxdt,
\endaligned
\end{equation*}
where $c_0$ is a fixed constant. The thesis follows by choosing $\tau=\frac{1}{2}$.
\end{proof}

\begin{prop}\label{qintegrabilitytime}
Let $u$ be a weak solution to \eqref{approx1}, then
\begin{equation*}
    \p_tu\in L^q_{loc}\left(Q\right),
\end{equation*}
for any $q\geq1$. Moreover, there is a constant $c=c(n,p,L)>0$ such that, for any $\beta\geq0$ and any non-negative $\eta\in C^1\left([0,T], C_c^\infty(\Omega)\right)$ vanishing on the parabolic boundary, it holds
\begin{equation}\label{Lqestimateptu}
    \int_0^T\int_\Om \eta^{\beta+2}|\p_t u|^{\beta+2}\ dxdt\leq c\ |\textnormal{supp}(\eta)|\left(M^{2p-2}\|\nabla_\e\eta\|_\infty^2+M^p\|\eta\p_t\eta\|_\infty\right)^{\frac{\beta+2}{2}},
\end{equation}
where $M=\sup_{\textnormal{supp}(\eta)}\left(\delta+|\nabla_\e u|^2\right)^{\frac{1}{2}}$.
\end{prop}

\begin{proof}
Since $u$ solves \eqref{approx1}, we can rewrite
\begin{equation}
    |\p_t u|^{\beta+2}=|\p_t u|^\beta\p_t u \sum_{i=1}^{2n+1}X^\e_i\left(D_if_\delta(\nabla_\e u)\right).
\end{equation}
From the previous identity and an integration by parts it follows
\begin{equation}\label{6marzo25_3}
\aligned
    L&:=\int_0^T\int_\Om\eta^{\beta+2}|\p_tu|^{\beta+2}\ dxdt=\sum_{i=1}^{2n+1}\int_0^T\int_\Om \eta^{\beta+2}|\p_tu|^{\beta}\p_tuX^\e_i\left(D_if_\delta(\nabla_\e u)\right)\ dxdt\\
    &=-(\beta+2)\sum_{i=1}^{2n+1}\int_0^T\int_\Om\eta^{\beta+1}|\p_tu|^{\beta}\p_tuD_if_\delta(\nabla_\e u)X^\e_i\eta\ dxdt\\
    &-(\beta+1)\sum_{i=1}^{2n+1}\int_0^T\int_\Om\eta^{\beta+2}|\p_tu|^{\beta}D_if_\delta(\nabla_\e u)X^\e_i\left(\p_t\eta\right)\ dxdt=I_1+I_2.
\endaligned
\end{equation}
To estimate them we use \eqref{structurecondition02} and the H\"older inequality. As for $I_1$, we obtain
\begin{equation}\label{7marzo25_5}
\aligned
    |I_1|&\leq(\beta+2)\int_0^T\int_\Om\eta^{\beta+1}|\p_tu|^{\beta+1}\left|\sum_{i=1}^{2n+1}\lambda_{i,\delta}(\nabla_\e u)X^\e_i u X_i\eta\right|\ dxdt\\
    &\leq (2n+1)(\beta+2)\int_0^T\int_\Om\eta^{\beta+1}|\p_tu|^{\beta+1}\left(\delta+|\nabla_\e u|^2\right)^{\frac{p-1}{2}}|\nabla_\e\eta|\ dxdt\\
    &\leq (2n+1)(\beta+2)\ |\textnormal{supp}(\eta)|^{\frac{1}{\beta+2}}\ \|\nabla_\e\eta\|_\infty\  M^{p-1}\ L^{\frac{\beta+1}{\beta+2}},
\endaligned
\end{equation}
where $M=\sup_{\textnormal{supp}(\eta)}\left(\delta+|\nabla_\e u|^2\right)^{\frac{1}{2}}$.

As for $I_2$, we obtain
\begin{equation}\label{7marzo25_2}
\aligned
    |I_2|&\leq (\beta+1)\int_0^T\int_\Om\eta^{\beta+2}|\p_tu|^{\beta}\left|\sum_{i=1}^{2n+1}\lambda_{i,\delta}(\nabla_\e u)X^\e_i u X_i\left(\p_tu\right)\right|\ dxdt\\
    &\leq (\beta+1)\int_0^T\int_\Om\eta^{\beta+2}|\p_tu|^{\beta}\left(\sum_{i=1}^{2n+1}\lambda_{i,\delta}(\nabla_\e u)(X^\e_i u)^2\right)^{\frac{1}{2}}\left(\sum_{i=1}^{2n+1}\lambda_{i,\delta}(\nabla_\e u)(X^\e_i\p_t u)^2\right)^{\frac{1}{2}}\ dxdt\\
    &(2n+1)^{\frac{1}{2}}(\beta+1)\int_0^T\int_\Om\eta^{\beta+2}|\p_tu|^{\beta}\left(\delta+|\nabla_\e u|^2\right)^{\frac{p}{4}}\left(\sum_{i=1}^{2n+1}\lambda_{i,\delta}(\nabla_\e u)(X^\e_i\p_t u)^2\right)^{\frac{1}{2}}\ dxdt\\
    &\leq (2n+1)^{\frac{1}{2}}(\beta+1)|\textnormal{supp}(\eta)|^{\frac{1}{\beta+2}}\ M^{\frac{p}{2}}\ L^{\frac{\beta}{2(\beta+2)}}\ J^{\frac{1}{2}},
\endaligned
\end{equation}
where 
\begin{equation*}
    J:=\int_0^T\int_\Om\eta^{\beta+4}|\p_tu|^\beta\sum_{i=1}^{n}\lambda_{i,\delta}(\nabla_\e u)\left((X^\e_i\p_t\eta)^2+(X^\e_{i+n}\p_t\eta)^2+(X^\e_{2n+1}\p_t\eta)^2\right)\ dxdt.
\end{equation*}
From \eqref{caccioppolidtu} it follows that 
\begin{equation}\label{7marzo25_3}
    J\leq c(M^{p-2}\|\nabla_\e\eta\|_\infty^2+\|\eta\p_t\eta\|_\infty)L,
\end{equation}
where $c=c(n,p,L)>0$. Hence, combining \eqref{7marzo25_2} and \eqref{7marzo25_3} we obtain
\begin{equation}\label{7marzo25_4}
    |I_2|\leq c\ (\beta+1)\ |\textnormal{supp}(\eta)|^{\frac{1}{\beta+2}}\ M^{\frac{p}{2}} \left(M^{p-2}\|\nabla_\e\eta\|_\infty^2+\|\eta\p_t\eta\|_\infty\right)^{\frac{1}{2}}L^{\frac{\beta+1}{\beta+2}},
\end{equation}
eventually changing the constant $c=c(n,p,L)>0$.

Combining \eqref{6marzo25_3}, \eqref{7marzo25_5} and \eqref{7marzo25_4} the thesis follows

\end{proof}

\section{Local Lipschitz regularity for solutions}\label{S:proofmainthm}

This section is devoted to the proof of the Main Theorem. The argument is based on passing to the limit in the estimates \eqref{mainestimateZu}, \eqref{moser-conclusion}, and \eqref{Lqestimateptu}, all of which remain stable as the approximating parameters $\varepsilon$ and $\delta$ tend to zero.

We begin by introducing sub-Riemannian parabolic cylinders, which are analogous to the Riemannian ones defined in Definition \ref{riemanniancylinder}. 

\begin{dfn}\label{cylinder}
Let $(x_0,t_0)\in \Om\times (0,T)$. For $r,\mu>0$, the parabolic cylinder $Q_{\mu,r}(x_0, t_0)\subset Q$ of center $(x_0,t_0)$ is the set 
$$Q_{\mu,r}(x_0, t_0)= B(x_0, r) \times (t_0-\mu\,r, t_0 ).$$

We call parabolic boundary of the cylinder $Q_{\mu,r}(x_0, t_0)\subset Q$
the set 
$$\p Q_{\mu,r}(x_0, t_0):=B(x_0, r) \times \{t_0-\mu\,r\} \cup \partial B(x_0, r) \times [t_0-\mu\,r, t_0).$$
\end{dfn}

To simplify the notation throughout the proof, we use a single approximation parameter $\varepsilon$ to denote both the sub-elliptic regularization and the Riemannian approximation.

\subsection{Proof of the Main Theorem}

The strategy for the proof is as follows: we begin by recalling some additional properties of the solution $u$; next, we approximate the equation and establish a uniform energy bound for the approximating solutions. Finally, we pass to the limit and transfer the desired estimates to the solution $u$.

\subsubsection*{Step 1: Additional properties of $u$.} Let us first recall that, for any $\tilde\Omega \subset\subset \Omega$ and any $0 < t_1 < t_2 < T$, we have
\begin{equation}\label{aggiuntive}
    \p_tu\in L^{p'}\left((t_1,t_2),HW^{-1,p'}(\tilde\Omega)\right)\quad\textnormal{and}\quad u\in C\left([t_1,t_2],L^2(\tilde\Omega)\right),
\end{equation}
where $p'=\frac{p}{p-1}$ and $HW^{-1,p'}(\Omega)$ denotes the dual space of $HW^{1,p}_0(\Omega)$. These are functional-analytic properties that do not depend on the geometry of the Heisenberg group $\mathbb{H}^n$, and thus hold as in the Euclidean setting. We briefly recall the argument for completeness.

For any $\tilde\Omega\subset\subset\Omega$ and $0<t_1<t_2<T$, using the weak formulation of equation \eqref{maineq01} and the fact that \(Df(\nabla_H u)\in L^{p'}_{\rm loc}(Q)\), we obtain
\[
\begin{split}
\left|\int_{t_1}^{t_2}\int_{\tilde\Omega} u\,\p_t\phi\,dxdt\right|&=\left|\int_{t_1}^{t_2}\int_{\tilde\Omega}\sum_{i=1}^{2n}D_if(\nabla_H u)X_i \phi\,dxdt\right|\\
&\le C\,\|\nabla_H \phi\|_{L^p(\tilde\Omega\times(t_1,t_2))}\leq \tilde{C} \|\phi\|_{L^p((t_1,t_2),HW^{1,p}_0(\tilde\Omega))}, \quad \forall\phi\in C^\infty_c(\tilde\Omega\times(t_1,t_2)).
\end{split}
\]
By density, the linear functional
\[
\Lambda:\phi\mapsto \int_{t_1}^{t_2}\int_{\tilde\Omega}u\,\p_t\phi\,dxdt,
\]
extends to the whole space \(L^p((t_1,t_2),HW^{1,p}_0(\tilde\Omega))\). This implies (see for instance \cite[Theorem 1.5, Chapter III]{Sh}) that
\[
\Lambda \in \Big(L^p((t_1,t_2),HW^{1,p}_0(\tilde\Omega))\Big)'=L^{p'}((t_1,t_2), HW^{-1,p'}(\tilde\Omega)).
\]
By the definitions of $\Lambda$ and of weak derivative, we obtain the first property in \eqref{aggiuntive}.

The second property in \eqref{aggiuntive} follows from the embedding result (see for instance \cite[Proposition 1.2, Chapter III]{Sh}):
\[
\Big\{\varphi\in L^p((t_1,t_2),HW^{1,p}_0(\tilde\Omega))\, :\, \p_t\varphi\in L^{p'}((t_1,t_2),HW^{-1,p'}(\tilde\Omega))\Big\}\subset C([t_1,t_2],L^2(\tilde\Omega)).
\]
In light of \eqref{aggiuntive}, and 
since we are only interested in local results, it is not restrictive to assume from the outset that
\begin{equation*}
    \p_tu\in L^{p'}\left((0,T),HW^{-1,p'}(\Omega)\right)\quad\textnormal{and}\quad u\in C\left([0,T],L^2(\Omega)\right).
\end{equation*}
\subsubsection*{Step 2: Regularization of the equation.}
By density, we consider a sequence of smooth functions $\left(\varphi_k\right)_{k\in\mathbb{N}}\subset C^\infty(Q)$ such that 
$$\varphi_k\to u \textnormal{ in }L^p\left((0,T),HW^{1,p}(\Omega)\right)\cap C\left([0,T],L^2(\Omega)\right),\quad\p_t\varphi_k\to\p_tu \textnormal{ in }L^{p'}\left((0,T),HW^{-1,p'}(\Omega)\right),$$
as $k\to\infty$. We then consider the following regularized Cauchy-Dirichlet problem: 
\begin{equation}\label{20aprile25}
\begin{cases}
    \p_t u^k_\e=\textnormal{div}_\e\left(Df_\e(\nabla_\e u^k_\e)\right),\quad &\textnormal{in }Q,\\
    u^k_\e=\varphi_k,\quad &\textnormal{in }\p Q.
\end{cases}
\end{equation}
The advantage of the Cauchy–Dirichlet problem \eqref{20aprile25} lies in the fact that it is Riemannian and therefore the existence of a unique solution $u^k_\e\in C^\infty(Q)\cap L^p\left((0,T),W^{1,p,\e}(\Omega)\right)$, such that \begin{equation}\label{24settembre2025}
    \p_t u^k_\e\in C^\infty(Q)\cap L^{p'}\left((0,T),W^{-1,p',\e}(\Omega)\right),
\end{equation}
where $W^{-1,p',\e}(\Omega)$ is the dual space of $W_0^{1,p,\e}(\Omega)$, is guaranteed by classical existence and regularity theory for Riemannian parabolic
equations; see for instance \cite{uraltsevaparabolic}.

\subsubsection*{Step 3: Uniform energy bound for approximating solutions.}
To obtain an energy estimate for $|\nabla_\varepsilon u^k_\varepsilon|^p$ that is uniform in $\varepsilon$ and $k$, we test the weak formulation of \eqref{20aprile25} with the function $\phi=\chi(t)\left(u^k_\e-\varphi_k\right)$, where $\chi\in C_c((0,T))$ satisfies $\chi\equiv1$ on $(0,t_0)$ and $\chi\equiv0$ on $(t_0,T)$. This yelds:
\begin{equation*}
    \int_0^{t_0}\int_\Omega u^k_\e\p_t\left(u^k_\e-\varphi_k \right)\ dxdt=\sum_{i=1}^{2n+1}\int_0^{t_0}\int_\Omega D_if_\e(\nabla_\e u^k_\e)X^\e_i \left(u^k_\e-\varphi_k\right)\ dxdt.
\end{equation*}
The left-hand side can be written as:
\begin{equation*}
    \aligned
    \int_0^{t_0}\int_\Omega u^k_\e\p_t\left(u^k_\e-\varphi_k \right)\ dxdt&=-\int_0^{t_0}\int_\Omega \p_t u^k_\e\left(u^k_\e-\varphi_k \right)\ dxdt\\&=-\frac{1}{2}\int_{\Omega}|u^k_\e-\varphi_k|^2(\cdot,t_0)\, dx-\int_0^{t_0}\int_\Omega\p_t\varphi_k\left(u^k_\e-\varphi_k \right)\, dxdt,
    \endaligned
\end{equation*}
where we integrated by parts in time. Therefore, we obtain: 
\begin{equation*}
    \aligned
    \frac{1}{2}\int_{\Omega}|u^k_\e-\varphi_k|^2(\cdot,t_0)\, dx&+\sum_{i=1}^{2n+1}\int_0^{t_0}\int_\Omega D_if_\e(\nabla_\e u^k_\e)X^\e_i \left(u^k_\e -\varphi_k\right)\ dxdt\\&=-\int_0^{t_0}\int_\Omega\p_t\varphi_k\left(u^k_\e-\varphi_k \right)\, dxdt.
    \endaligned
\end{equation*}
By Lemma \ref{8agosto2025}, for each $i=1,\ldots,2n+1$, we have:
\begin{equation*}
    \aligned
    &\frac{1}{2}\int_{\Omega}|u^k_\e-\varphi_k|^2(\cdot,t_0)\, dx+\int_0^{t_0}\int_\Omega|X^\e_i u^k_\e|^p\ dxdt\\&\leq c\bigg{(}\frac{1}{2}\int_{\Omega}|u^k_\e-\varphi_k|^2(\cdot,t_0)\, dx+\sum_{i=1}^{2n+1}\int_0^{t_0}\int_\Omega D_if_\e(\nabla_\e u^k_\e)X^\e_i \left(u^k_\e -\varphi_k\right)\ dxdt\\&+\int_0^{t_0}\int_\Omega(\e+|\nabla_\e\varphi_k|^2)^{\frac{p}{2}}\ dxdt\bigg{)}\\
    &=c\left(-\int_0^{t_0}\int_\Omega\p_t\varphi_k\left(u^k_\e-\varphi_k \right)\, dxdt+\int_0^{t_0}\int_\Omega(\e+|\nabla_\e\varphi_k|^2)^{\frac{p}{2}}\ dxdt\right) \\&\leq c\left(\tau\int_0^{t_0}\int_\Omega|\nabla_\e u^k_\e-\nabla_\e\varphi_k|^p\, dxdt +c_\tau\|\p_t\varphi_k\|_{L^{p'}\left((0,t_0),HW^{-1,p'}(\Omega)\right)}+\int_0^{t_0}\int_\Omega(\e+|\nabla_\e\varphi_k|^2)^{\frac{p}{2}}\ dxdt\right)\\
    &\leq c\tau\int_0^{t_0}\int_\Omega|\nabla_\e u^k_\e|^p\, dxdt + c_\tau\|\p_t\varphi_k\|_{L^{p'}\left((0,t_0),HW^{-1,p'}(\Omega)\right)}+c\int_0^{t_0}\int_\Omega(\e+|\nabla_H\varphi_k|^2)^{\frac{p}{2}}\ dxdt+o(\e)\\&\leq c\tau\int_0^{t_0}\int_\Omega|\nabla_\e u^k_\e|^p\, dxdt + c_\tau\|\p_tu\|_{L^{p'}\left((0,t_0),HW^{-1,p'}(\Omega)\right)}\\&+c\int_0^{t_0}\int_\Omega(\e+|\nabla_Hu|^2)^{\frac{p}{2}}\ dxdt+o(\e)+o\left(\frac{1}{k}\right),
    \endaligned
\end{equation*}
where $\tau>0$ and  $c,c_\tau>0$ depend only on $p$. Summing over $i=1,\ldots,2n+1$ and choosing $\tau$ small enough, we conclude that the following inequality holds true for almost every $t_0\in(0,T)$:
\begin{equation}\label{2maggio25}
    \aligned
    \int_{\Omega}|u^k_\e-\varphi_k|^2(\cdot,t_0)\, dx&+\int_0^{t_0}\int_\Omega|\nabla_\e u^k_\e|^p\, dxdt\\&\leq c\int_0^{t_0}\int_\Omega\left(\e+|\nabla_H u|^2\right)^\frac{p}{2}\, dxdt+c\|\p_tu\|_{L^{p'}\left((0,t_0),HW^{-1,p'}(\Omega)\right)}+o(\e)+o\left(\frac{1}{k}\right),
    \endaligned
\end{equation}
where $c=c(n,p)>0$ is independent of $\e$ and $k$. 

\subsubsection*{Step 4: Convergence of the approximation scheme}
The Riemannian Poincaré inequality, together with the uniform energy bound \eqref{2maggio25}, implies the existence of a constant $M_0 > 0$, independent of $\varepsilon$ and $k$, such that
\begin{equation*}
    \|u_\varepsilon^k\|_{L^p((0,T), W^{1,p,\varepsilon}(\Omega))} \leq M_0.
\end{equation*}
Therefore, the sequence $\left(u^k_\varepsilon\right)_{\varepsilon,k}$ is bounded in $L^p((0,T), W^{1,p,\varepsilon}(\Omega))$.

On the other hand, Theorem \ref{moser} implies that $\nabla_\varepsilon u_\varepsilon^k \in L^\infty_{\mathrm{loc}}(Q, \mathbb{R}^{2n+1})$. Moreover, the uniform energy bound \eqref{2maggio25} ensures the existence of a constant $M_1 > 0$, independent of $\varepsilon$ and $k$, such that for any $\tilde\Omega \subset\subset \Omega$ and $0 < t_1 < t_2 < T$,
\begin{equation}\label{14agosto}
    \|\nabla_\varepsilon u^k_\varepsilon\|_{L^\infty(\tilde\Omega \times (t_1,t_2), \mathbb{R}^{2n+1})} 
    \leq c\mu^{\frac{1}{2}}\, \max\left\{ \left( \frac{1}{\mu r^{N+2}} \int_0^T \int_\Omega (\varepsilon + |\nabla_\varepsilon u^k_\varepsilon|^2)^{\frac{p}{2}}\, dxdt \right)^{\frac{1}{2}}, \mu^{\frac{p}{2(2-p)}} \right\} \leq M_1.
\end{equation}
Hence, the sequence $\left(\nabla_\varepsilon u^k_\varepsilon\right)_{\varepsilon,k}$ is bounded in $L^\infty_{\mathrm{loc}}(Q, \mathbb{R}^{2n+1})$.

In addition, Proposition \ref{Poincaretypeineq} implies that $Z u^k_\varepsilon \in L^p_{\mathrm{loc}}(Q)$. Furthermore, the uniform energy estimate \eqref{2maggio25} yields the existence of a constant $M_2 > 0$, independent of $\varepsilon$ and $k$, such that for any $\tilde\Omega \subset\subset \Omega$ and $0 < t_1 < t_2 < T$,
\begin{equation*}
    \|Z u^k_\varepsilon\|_{L^p(\tilde\Omega \times (t_1,t_2))} \leq M_2.
\end{equation*}
Therefore,
\begin{equation}\label{14agosto1}
    \|X_{2n+1}^\varepsilon u^k_\varepsilon\|_{L^p(\tilde\Omega \times (t_1,t_2))} = \varepsilon \|Z u^k_\varepsilon\|_{L^p(\tilde\Omega \times (t_1,t_2))} \to 0 \quad \text{as } \varepsilon \to 0^+.
\end{equation}

Hence, we can conclude that, up to a subsequence, there exists a function $u_0 \in L^p((0,T), HW^{1,p}(\Omega))$, with $\nabla_H u_0 \in L^\infty_{\mathrm{loc}}(Q, \mathbb{R}^{2n})$, such that for any $\tilde\Omega \subset\subset \Omega$ and $0 < t_1 < t_2 < T$,
\begin{equation*}
    u^k_\varepsilon \rightharpoonup u_0 \quad \text{in } L^p((0,T), HW^{1,p}(\Omega)), \quad \text{and} \quad \nabla_\varepsilon u^k_\varepsilon \rightharpoonup^* \nabla_H u_0 \quad \text{in } L^\infty(\tilde\Omega \times (t_1,t_2), \mathbb{R}^{2n+1}),
\end{equation*}
as $\varepsilon \to 0^+$ and $k \to \infty$.

Moreover, the weak convergence in $L^p((0,T), HW^{1,p}(\Omega))$, together with the fact that $u^k_\varepsilon - \varphi_k \in L^p((0,T), HW^{1,p}_0(\Omega))$, implies that $u_0 - u \in L^p((0,T), HW^{1,p}_0(\Omega))$, and hence $u = u_0$ on the lateral boundary of $Q$. From \eqref{2maggio25}, it also follows that $u(\cdot,0) = u_0(\cdot,0)$, since
\begin{equation*}
\aligned
    \int_{\Omega}|u_0-u|^2(\cdot,t_0)\, dx&\leq \liminf_{\substack{\e\to0^+\\k\to\infty}}\int_{\Omega}|u^k_\e-\varphi_k|^2(\cdot,t_0)\, dx\\&\leq c\int_0^{t_0}\int_\Omega\left(\e+|\nabla_H u|^2\right)^\frac{p}{2}\, dxdt+c\|\p_tu\|_{L^{p'}\left((0,t_0),HW^{-1,p'}(\Omega)\right)},
\endaligned
\end{equation*}
and the claim follows by letting $t_0$ tend to $0^+$.

From \eqref{standardcaccioppoli2} and \eqref{14agosto}, it follows that 
\[
\nabla_\varepsilon\left(|X^\varepsilon_\ell u_\varepsilon^k|^{\beta+p}\right) \in L^2_{\mathrm{loc}}(Q, \mathbb{R}^{2n+1}), \quad \text{for all } \ell = 1, \ldots, 2n+1.
\]
In particular, there exists a constant $M_3 > 0$, independent of $\varepsilon$ and $k$, such that for any $\tilde\Omega \subset\subset \Omega$ and $0 < t_1 < t_2 < T$,
\[
\left\|\nabla_\varepsilon\left(|X^\varepsilon_\ell u_\varepsilon^k|^{\beta+p}\right)\right\|_{L^2(\tilde\Omega \times (t_1, t_2))} \leq M_3, \quad \forall \ell = 1, \ldots, 2n+1.
\]

By virtue of Theorem \ref{Rellichthm} and \eqref{24settembre2025}, we can apply \cite[Theorem 5]{Simon} with $X = W^{1,p,\varepsilon}(\tilde\Omega)$, $B = L^r(\tilde\Omega)$ for any $r < p^* = \frac{Np}{N - p}$, and $Y = W^{-1,p',\varepsilon}(\tilde\Omega)$. We conclude that there is $v\in L^r_{\mathrm{loc}}(Q)$ such that
\[
|X^\varepsilon_\ell u_\varepsilon^k|^{\beta+p} \to v \quad \text{in } L^r_{\mathrm{loc}}(Q), \quad \text{for any } r < p^* \text{ and any } \beta \geq 0.
\]
Consequently, we deduce that there is $\tilde v\in L^q_{loc}(\Omega)$ such that
\[
X^\varepsilon_\ell u_\varepsilon^k \to \tilde v \quad \text{in } L^q_{\mathrm{loc}}(Q), \quad \text{for any } 1 \leq q < \infty.
\]

From the previous discussion and \eqref{14agosto1}, we conclude that
\[
\nabla_\varepsilon u_\varepsilon^k \to \nabla_H u_0 \quad \text{in } L^p_{\mathrm{loc}}(Q, \mathbb{R}^{2n+1}),
\]
and hence, up to a subsequence,
\[
\nabla_\varepsilon u_\varepsilon^k \to \nabla_H u_0 \quad \text{a.e. in any compact subset } K \subset\subset \Omega.
\]

Our next goal is to prove that $u_0$ is a weak solution to equation \eqref{maineq}. If this is the case, then since $u_0 = u$ on $\partial Q$, it follows by the comparison principle that $u = u_0$ in $Q$. To this end, we pass to the limit in the weak formulation of equation \eqref{approx1}:
\begin{equation*}
    \int_0^T \int_\Omega u^k_\e\, \p_t\phi\  dx dt=\int_0^T\int_{\Omega}\sum_{i=1}^{2n+1}D_i f_{\e}(\nabla_\e u^k_\e)X_i^\e\phi\ dxdt,\quad \forall \phi\in C^\infty_c(Q).
\end{equation*}
This passage to the limit is straightforward: from \eqref{14agosto}, $|\nabla_\e u_\e^k(x,t)|$ is uniformly bounded for $(x,t)\in\textnormal{supp }\phi$, and hence $|Df_\e(\nabla_\e u_\e^k(x,t))|$ is uniformly bounded as well. Moreover, the almost everywhere convergence $\nabla_\e u^k_\e(x,t)\to \nabla_H u_0(x,t)$ in $\textnormal{supp }\phi$, together with the convergence $Df_\e(z)\to Df(z_H)$, implies that \[Df_\e(\nabla_\e u^k_\e(x,t))\to Df(\nabla_H u_0(x,t))\quad \text{a.e. in } \text{supp } \phi.\]
We can then apply the dominated convergence theorem to the second integral and use the weak convergence of $u_\e^k$ to $u_0$ in $L^p((0,T),HW^{1,p}(\Omega))$ to conclude:
\begin{equation*}
    \int_0^T \int_\Omega u_0\,\p_t\phi\  dx dt=\int_0^T\int_{\Omega}\sum_{i=1}^{2n}D_i f(\nabla_H u_0)X_i\phi\ dxdt,\quad \forall \phi\in C^\infty_c(Q).
\end{equation*}

\subsubsection*{Step 5: Local Lipschitz estimate for $u$.}
Theorem \ref{moser}, combined with the almost everywhere convergence $\nabla_\varepsilon u^k_\varepsilon(x,t) \to \nabla_H u_0(x,t)$ and the convergence of Riemannian balls to sub-Riemannian ones in the Hausdorff sense, implies that 
\begin{equation*}
    \aligned
    \|\nabla_Hu\|_{L^\infty(Q_{\mu,r})}&\leq\liminf_{\e\to0^+}\|\nabla_\e u^k_\e\|_{L^\infty(Q^\e_{\mu,r},\mathbb{R}^{2n+1})}\\&\leq c\mu^{\frac{1}{2}}\,\liminf_{\e\to0^+}\left\{\max\left\{\left(\frac{1}{\mu r^{N+2}}\int\int_{Q^\e_{\mu,2r}}\left(\e+|\nabla_\e u^k_\e|^2\right)^{\frac{p}{2}}\ dxdt\right)^{\frac{1}{2}},\mu^{\frac{p}{2(2-p)}}\right\}\right\}\\
    &=c\mu^{\frac{1}{2}}\,\max\left\{\left(\frac{1}{\mu r^{N+2}}\int\int_{Q_{\mu,2r}}|\nabla_H u|^p\ dxdt\right)^{\frac{1}{2}},\mu^{\frac{p}{2(2-p)}}\right\}.
    \endaligned
\end{equation*}

\subsubsection*{Step 6: Local $L^q$ estimate for $Zu$ and $\p_tu$.}
Finally, Proposition \ref{Poincaretypeineq}, together with \eqref{14agosto}, implies that the sequence $\left(Z u^k_\varepsilon\right)_{\varepsilon,k}$ is uniformly bounded in $L^q_{\text{loc}}(Q)$ for all $1 \leq q < \infty$. Therefore, up to a subsequence, for any $\tilde\Omega \subset\subset \Omega$ and $0 < t_1 < t_2 < T$, there exists $v \in L^q_{\text{loc}}(Q)$ such that $Z u^k_\varepsilon \rightharpoonup v$ in $L^q(\tilde\Omega \times (t_1,t_2))$. 

From the definition of weak derivative and the convergence results established earlier, we conclude - possibly after passing to a further subsequence - that $v$ coincides with the vertical derivative of the limit function, i.e., $v = Z u$. 

An analogous argument applies to the sequence $\left(\partial_t u^k_\varepsilon\right)_{\varepsilon,k}$, for which the local $L^q$ uniform bound \eqref{Lqestimateptu} holds. Hence, up to a subsequence, we also obtain that $\partial_t u^k_\varepsilon \rightharpoonup \partial_t u$ in $L^q_{\text{loc}}(Q)$ for all $1 \leq q < \infty$.

\printbibliography

\end{document}